\documentclass[reqno, 10pt]{amsart}  
\usepackage{amsmath}
\usepackage{amssymb}
\usepackage{hyperref}
\hypersetup{colorlinks=true}
\usepackage{graphicx}
\graphicspath{~/Desktop/first-order hyperbolic/}
\numberwithin{equation}{section}

\newtheorem{prop}{Proposition}[section]
\newtheorem{lemma}[prop]{Lemma}
\newtheorem{thm}[prop]{Theorem}
\newtheorem{cor}[prop]{Corollary}

\newcommand{\shrink}[1]{ {\scriptstyle {\textstyle {#1} } } }
\newcommand{\smfrac}[2]{ \shrink{ \frac{#1}{#2} } }

\newcommand{\dom}{\mathrm{dom}} 

\newcommand{\blin}{\boldsymbol{\big \langle}} 
\newcommand{\brin}{\boldsymbol{\big \rangle}}

\newcommand{\lin}{\langle}
\newcommand{\rin}{\rangle}
\newcommand{\Sym}{\mathbb{S}^{n}}
\newcommand{\Symp}{\Sym_{\plus}} 
 
\newcommand{\realsnp}{\mathbb{R}^n_{\plus}}

\newcommand{\bdy}{\mathrm{bdy}}

\newcommand{\plus}{{\scriptscriptstyle +}}

\renewcommand{\int}{\mathrm{int}}

\newcommand{\epi}{\textrm{epi}}

\newcommand{\cp}{\mathrm{CP}} 

\newcommand{\lambdamin}{\lambda_{\mathrm{min}}}

\newcommand{\diamval}{\mathrm{diam}_{z}  }

\usepackage{lipsum}

\begin{document}

\newpage 
$ \textrm{~} $ \quad \vspace{-3mm}

\title[``Efficient'' General Subgradient-Methods]{``Efficient'' Subgradient Methods
\\ for General Convex Optimization}

\begin{abstract}
A subgradient method is presented for solving general convex optimization problems, the main requirement being that a strictly-feasible point is known.  A feasible sequence of iterates is generated, which converges to within user-specified error of optimality.  Feasibility is maintained with a line-search at each iteration, avoiding the need for orthogonal projections onto the feasible region (the operation that limits practicality of traditional subgradient methods).  Lipschitz continuity is not required, yet the algorithm is shown to possess a convergence rate analogous to rates for traditional methods, albeit with error measured relatively, whereas traditionally error has been absolute.  The algorithm is derived using an elementary framework that can be utilized to design other such algorithms.
\end{abstract} \vspace{-3mm}

\author[J. Renegar]{James Renegar  }
\address{School of Operations Research and Information Engineering,
 Cornell University, Ithaca, NY, U.S.}
 
\thanks{Research supported in part by NSF CCF Award \#1552518. }

\maketitle

\vspace{-7mm}

\section{{\bf  Introduction}}  \label{sect.a}
Given a convex optimization problem possessing little more structure than having a known strictly-feasible point, we provide a  simple transformation to an equivalent convex optimization problem which has only linear equations as constraints, and has Lipschitz-continuous objective function defined on the whole space. Virtually any subgradient method can be applied to solve the equivalent problem.  Relying only on a line-search during an iteration, the resulting iterate  is made to be a feasible point for the original problem. Moreover, this sequence of feasible points converges to within user-specified error of optimality.   \vspace{1mm}

The algorithm can be implemented directly in terms of the original problem, leaving the equivalent problem hidden from view. 
In this introductory section we present two examples of algorithms and complexity bounds that ensue when standard subgradient methods are applied to the equivalent problem. Not until subsequent  sections is the equivalent problem revealed. \vspace{1mm}

Let $ {\mathcal E} $ denote a finite-dimensional real vector space with inner product $ \langle \; , \; \rangle $ and associated norm $ \| \, \, \| $. \vspace{1mm}

Consider an optimization problem
\begin{equation}   \label{eqn.aa}
  \begin{array}{rl} \min   & f(x) \\
                           \textrm{s.t.} & x \in \mathrm{Feas} \; . \end{array} 
                           \end{equation}  
Assume $ f: {\mathcal E} \rightarrow  (- \infty, \infty] $ is an extended-valued, lower-semicontinous convex function.  Equivalently, assume $ f: {\mathcal E} \rightarrow \mathbb{R} $ has closed and convex epigraph, $ \mathrm{epi}(f) := \{ (x,t) \mid x \in \dom(f), \,  t \geq f(x) \} $,  (where $ \dom(f) $ is the (effective) domain of $ f $ (i.e., the set on which $ f $ is finite)).  \vspace{1mm}

Let
\[ \mathrm{Feas} = \{ x \in S \mid Ax = b \} \; , \]
where $ S $ is a closed, convex set with nonempty interior.  \vspace{1mm}

 Assume $ \bar{e} $ is a known feasible point contained in $ \mathrm{int}( S \cap \mathrm{\mathrm{dom}(f)}) $.   \vspace{1mm}

 Let $ \bar{D}  $ denote the diameter of the sublevel set $ \{ x \in \mathrm{Feas}: f(x) \leq f( \bar{e} ) \} $ (i.e.,  the supremum of distances between pairs of points in the set).  We assume $ \bar{D}  $ is finite, which together with previous assumptions implies the optimal value of (\ref{eqn.aa})  is attained at some feasible point.  Let $ f^* $ denote the optimal value. \vspace{1mm}

Let $ \hat{f} $ be a user-chosen constant satisfying $ \hat{f} > f( \bar{e} ) $ (thus, $ (\bar{e}, \hat{f}) \in \mathrm{int}(\mathrm{epi}(f)) $).  Define
\[ 
  \hat{r}  := \sup \{  r \mid \| x - \bar{e}  \| \leq r \textrm{ and }   Ax = b \,  \, \Rightarrow  \, \, x \in S \textrm{ and }   f(x) \leq \hat{f} \}  \; . \] 
Note $ \hat{r}  > 0 $, because $ \bar{e}  \in \mathrm{int}\big( S \cap  \mathrm{dom}(f)\big) $.   \vspace{1mm}

The scalars $ \bar{D}  $ and $ \hat{r}  $ appear in the complexity bounds for the representative algorithms presented below, but their values are not assumed to be known.  We now introduce notation to be used in specifying the algorithms.   \vspace{1mm}                

For fixed $ x \in {\mathcal E} $ and $ t \in \mathbb{R} $, and for scalars $ \alpha $, let
\[ x(\alpha) :=  \bar{e}  + \alpha \cdot (x - \bar{e} ) \quad \textrm{and} \quad t( \alpha) :=  \hat{f} + \alpha \cdot ( t - \hat{f}) \; .     \]
Define  
\[   \alpha_1(x) := \sup \{ \alpha \mid x(\alpha) \in S \} \; , \]
and
\[  \alpha_2(x,t) := \sup \{ \alpha \mid  f(x( \alpha)) \leq  t( \alpha) \} \; . \]
Note that $ \alpha_1(x), \alpha_2(x,t) > 0 $, due to our requirements for $ \bar{e}  $ and $ \hat{f} $.   \vspace{1mm}

We are concerned only with  pairs $ (x,t) $ for which $ Ax = b $ and $ t < \hat{f} $.  (Observe then, $ A x( \alpha) = b $ and $ t(\alpha) < \hat{f} $ for all $ \alpha \geq 0 $.)  We claim in this setting, at least one of the values $ \alpha_1(x) $, $ \alpha_2(x,t) $ is finite. To establish the claim, assume $ \alpha_2(x,t) = \infty $ -- we show it easily follows that $ \alpha_1(x) $ is finite. \vspace{1mm}

Since $ \alpha_2(x,t) = \infty $, the convex univariate function
\begin{equation}  \label{eqn.ab} 
   \phi(\alpha) := f( x( \alpha)) - t( \alpha) 
   \end{equation}
has negative value for all $ \alpha \geq 0 $.  Hence, since $ t < \hat{f} $ (by assumption) -- and thus $ t( \alpha) \rightarrow - \infty $ as $ \alpha \rightarrow \infty $ --  the function $ \alpha \mapsto f( x( \alpha)) $ is unbounded below on the interval $ [0, \infty ) $. In particular, $ f(x( \alpha)) < f^* $ for some $ \alpha > 0 $, implying for this value of $ \alpha $,  $ x( \alpha) \notin \mathrm{Feas} $, that is, $ \alpha_1(x) < \alpha $, establishing the claim.  \vspace{1mm}

Let
\[  \alpha(x,t) := \min \{ \alpha_1(x), \alpha_2(x,t) \} \; , \]
a value we now know to be positive and finite (assuming $ Ax = b $ and $ t < \hat{f} $). For given pairs $ (x,t) $, the algorithms require that the scalar $ \alpha(x,t) $ be computed.  (This is the line search referred to above.)  In general, of course, $ \alpha(x,t) $  cannot be computed exactly.  However, there is hope for quickly gaining a ``good enough'' approximation to $ \alpha(x,t) $. \vspace{1mm}

To understand, first observe for many optimization problems, for each $ x $ the value $ \alpha_1(x) $ can readily be accurately approximated, or determined to be $ \infty $.  All that is required here is finding whether a given half-line, with endpoint in $ \mathrm{int}(S) $,   intersects the boundary of $ S $, and if so, accurately approximating where the intersection occurs.  For most feasible regions, this is far more easily accomplished than one of the main operations underlying traditional subgradient methods, where an (arbitrary) point outside $ \mathrm{Feas} $ must be orthogonally projected onto $ \mathrm{Feas} $. \vspace{1mm}

Assume for the moment that $ \alpha_1(x) $ has been accurately approximated or determined to equal $ \infty $.  For definiteness, assume $ \alpha_1(x) = \infty $, in which case we know $ \alpha_2(x,t) $ is finite.  Approximating $ \alpha_2(x,t) $ is the same as approximating the solution to $ \phi(\alpha) = 0 $, where $ \phi $ is the function (\ref{eqn.ab}).  As $ \phi $ is a convex function satisfying $ \phi (0) < 0 $, approximating the root can be accomplished by first finding a value $ \alpha $ sufficiently large so as to satisfy $ \phi(\alpha) > 0 $, and then proceeding to isolate the root using bisection.  \vspace{1mm}

On the other hand, if $ \alpha_1(x) $ is finite, then to determine whether $ \alpha_2(x,t) $ even needs to be computed, check whether $ \phi( \alpha_1(x)) \leq  0 $.  If so, there is no need to compute $ \alpha_2(x,t) $, because $ \alpha(x,t) = \alpha_1(x) $.  If instead, $ \phi(\alpha_1(x)) > 0 $, then to compute $ \alpha_2(x,t) $ proceed by bisection, with initial \vspace{1mm} interval $ [0, \alpha_1(x)] $. \vspace{1mm}

That bisection can be applied gives reason to hope a good approximation to $ \alpha(x,t) $ can be quickly computed.  Nonetheless, questions abound regarding what constitutes a ``good enough'' approximation.  In what follows, we duck the issue, simply assuming $ \alpha(x,t) $ can be computed exactly. \vspace{1mm}

For $ x,t $ satisfying $ Ax = b $ and $ t < \hat{f} $, let 
\begin{equation}  \label{eqn.ac} 
   \pi'(x,t) :=  (\bar{e} , \hat{f}  ) + \alpha(x,t) \cdot ( (x,t) - (\bar{e} , \hat{f})) \; . 
   \end{equation} 
Thus, in traveling from $ (\bar{e} , \hat{f}) $ in the direction $  (x,t) - (\bar{e} , \hat{f}) $, $ \pi'(x,t) $ is the first point $ (x',t') $ encountered for which either $ x' \in \mathrm{bdy}(S)  $ (boundary) or $ (x',t') \in \bdy(\epi(f))  $.  Note in particular that if $ (x',t') = \pi' (x,t) $, then $ x' $ is both feasible and lies in the domain of $ f $.
\vspace{1mm}

For $ x' \in \mathrm{bdy}(S) $, let
\[ G_1(x') := \left\{ \frac{-1}{ \lin v , \bar{e}  - x' \rin } \, v \mid \,  \vec{0} \neq  v \in N_S(x') \right\} \; , 
\]
where $ N_S(x') $ is the normal cone to $ S $ at $ x' $ ($\{ v \mid \forall x \in S, \lin v, x - x' \rin \leq 0 \} $).  The denominator is negative, because $ \bar{e}  \in \mathrm{int}(S) $.  \vspace{1mm}  

For $ (x',t') \in \mathrm{\bdy}(\mathrm{epi}(f)) $, let         
\[ G_2(x',t') := \left\{ \frac{-1}{\lin v, \bar{e} - x' \rin + ( \hat{f} - t') \delta  } \, v \mid  \, (\vec{0},0) \neq (v,\delta) \in N_{\mathrm{epi}(f)}(x',t')   \right\} \; , \] 
where here, normality is with respect to the inner product that assigns pairs $ (x_1, t_1) $, $ (x_2, t_2) $ the value $ \lin x_1, x_2 \rin + t_1 t_2 $. \vspace{1mm}

If $ (x',t') \in \mathrm{bdy}(\mathrm{dom}(f)) $ and $ x' \in \mathrm{int}(\mathrm{dom}(f)) $ -- hence $ t' = f(x') $ -- a more concrete description applies:
\begin{align*} 
  & x' \in \mathrm{int}(\mathrm{dom}(f)) \,  \\
& \qquad   \Rightarrow \, G_2(x', f(x')) =   
\left\{ \frac{1}{\hat{f} - ( f(x') + \lin v, \bar{e}  - x' \rin )}\, v \mid \,    v \in \partial f(x') \right\} \; ,  
\end{align*}  
where $ \partial f(x') $ is the subdifferential of $ f $ at $ x' $ (set of subgradients).  We  rely only on sets $ G_2(x',t') $ for which $ x' \in \mathrm{Feas} $ -- thus,   if  $ \mathrm{Feas} \cap \mathrm{dom}(f) \subseteq \mathrm{int}(\mathrm{dom}(f)) $, the more concrete description  applies.  (In the traditional literature on subgradient methods, the stronger condition  $ \mathrm{Feas} \subset \mathrm{int}(\mathrm{dom}(f)) $ is generally assumed, and hence the more concrete description certainly applies.) \vspace{1mm}

Finally, letting $ \bdy_1 := \{ (x',t') \mid x' \in \bdy(S) \} $ and $ \bdy_2 := \bdy(\epi(f)) $, define for $ (x',t') \in \bdy_1 \cup \bdy_2 \; , $ 
\[ G(x',t') := \begin{cases}  G_1(x')  & \textrm{if $ (x',t') \in \bdy_1 \setminus \bdy_2 \; , $}  \\
G_2(x',t') & \textrm{if $ (x',t') \in \bdy_2 \setminus \bdy_1\; , $} \\
\textrm{hull}(G_1(x') \, \cup \, G_2(x',t') ) & \textrm{if $ (x',t') \in \bdy_1 \cap  \bdy_2 \; , $}   
\end{cases}    \]
where ``hull'' denotes the convex hull. \vspace{1mm}

Let $ P $ denote orthogonal projection onto  the kernel of $ A $.   \vspace{1mm}

Following is a representative algorithm that arises from the framework developed in subsequent sections.  

\noindent 
\hrulefill
 
\noindent {\bf Algorithm A} \vspace{1mm}

\noindent 
\begin{tabular}{lll}
    (0) & Input: & $ 0 < \epsilon < 1 $ ,  \\
&  & $ \bar{e}  \in \mathrm{int}(S \cap \mathrm{dom}(f)) $ satisfying $ A \bar{e}  = b      $, and \\
&&   $ \hat{f}  $, a scalar satisfying $ \hat{f} > f( \bar{e} ) $.    \\
& Initialize: & $ (x_0, t_0) := (\bar{e} , f( \bar{e} )) $ and $ (x_0',t_0') :=  (\bar{e} , f( \bar{e} )) $. \\
(1) & Iterate: & Compute $ \tilde{x}_{k+1}  := x_k - \smfrac{\epsilon }{2 \| P(g') \|^2} P(g') $, \,   where $ g' \in G( x_k', t_k') $. \\
&& Let $ \alpha_{k+1} := \alpha( \tilde{x}_{k+1}, t_k) $ and $ (x_{k+1}', t_{k+1}') := \pi' (\tilde{x}_{k+1} , t_k) $. \\ 
&& If $ \alpha_{k+1} \geq 4/3 $, then let $ (x_{k+1}, t_{k+1}) :=  (x_{k+1}', t_{k+1}') $, \\
&& $ \textrm{~} $ \quad else let  $ (x_{k+1}, t_{k+1})  :=  (\tilde{x}_{k+1} , t_k) $.   \vspace{1mm} 
\end{tabular}  

\noindent 
\hrulefill
\vspace{2mm}

Critical to understanding the algorithm is the identity $ (x_k', t_k') = \pi'(x_k,t_k) $.  Iterates in the sequence $ \{ x_k \} $ need not be feasible nor in the domain of $ f $, whereas $ x_k' $ is feasible and in the domain.  \vspace{1mm}

Unlike traditional subgradient methods, no orthogonal projections onto $ \mathrm{Feas}  $ are required.
The only projections are onto the kernel of $ A $, projections which are computed efficiently if the range space is low dimensional, and if as a preprocessing step the linear operator $ (A A^*)^{-1} $ is computed and stored in memory (where $ A^* $ is the adjoint of $ A $, and where we are assuming, of course, that $ A $ is surjective, for the inverse to exist). \vspace{1mm}  

\begin{thm}  \label{thm.aa} 
For the feasible sequence $ \{ x_k' \} $ generated by Algorithm A,
\begin{align}  
 & \ell  \geq 8 \, \left(  \frac{\bar{D} }{\hat{r} } \right)^2    \, \left( \, \frac{1}{\epsilon^2} \, + \, \frac{1}{\epsilon} \, \log_{4/3} \left( 1 + \frac{\bar{D} }{\hat{r} } \right)       \, \right)     \label{eqn.ad}  \\
&  \quad                \quad \Rightarrow \quad \min_{k \leq \ell } \,   \frac{f(x_k') - f^*}{\hat{f} - f^*} \, \leq \, \epsilon \; . \label{eqn.ae} 
 \end{align}
 \end{thm} 
\vspace{2mm}

This result differs markedly from the traditional literature, in that there appears no Lipschitz constant for $ f $.  The theorem applies even when for no neighborhood of an optimal solution does a (local) Lipschitz constant exist, such as occurs, for example, when $ \mathrm{Feas} = \mathbb{R}^{2} $ and
\[ f(x_1,x_2) = \begin{cases} 0 & \textrm{if $ x = \vec{0} $}, \\
  x_1^2 + x_2^2/x_1 & \textrm{if $ x_1 > 0 $}, \\
  \infty & \textrm{otherwise}.    
\end{cases}  \]
(For every $ r > 0 $, the function fails to be Lipschitz continuous on the relatively-open set $ \{ x \in \dom(f) \mid \| x \| < r \} $.)
\vspace{1mm}

The equivalent problem is hidden from view, the problem to which a standard subgradient method is applied and then translated as Algorithm A.  The objective function for the equivalent problem  has Lipschitz constant bounded above by $ 1/ \hat{r} $. Thus, a Lipschitz constant does in fact appear in (\ref{eqn.ad}), even though $ f $ need not be Lipschitz continuous. \vspace{1mm}   
   
A tradeoff to gaining independence from Lipschitz continuity is that the error (\ref{eqn.ae})  is measured relatively rather than absolutely.  Dependence on $ 1/ \epsilon^2 $ occurs both in (\ref{eqn.ad})  and in the \vspace{1mm} traditional literature. \vspace{1mm}

 The criterion ``$ \alpha_{k+1} \geq 4/3 $'' is satisfied at most $ \log_{4/3}(1 + \bar{D}/ \hat{r}) $ times, leading to the logarithmic term in (\ref{eqn.ad}). The primary importance of using the criterion to create two cases is due to the accuracy needed in solving the equivalent problem being dependent not only on $ \epsilon $, but also on $ f^* $.  Use of the criterion provides a means to compensate for not knowing $ f^* $. \vspace{1mm}

When $ f^* $ is known, a streamlined algorithm can be devised. This algorithm does not require $ \epsilon $ as input.
\vspace{1mm}

\noindent 
\hrulefill
 
\noindent {\bf  Algorithm B} \vspace{1mm}

\noindent 
\begin{tabular}{lll}
    (0) & Input: & $ f^* $, the optimal value, \\
    & &   $ \bar{e}  \in \mathrm{int}( S \cap \mathrm{dom}(f)) $ satisfying $ A \bar{e}  = b      $, and \\
&&   $ \hat{f}  $, a scalar satisfying $ \hat{f} > f( \bar{e} ) $.    \\
& Initialize: & $ x_0 := \bar{e} $, $ (x_0',t_0') :=  (\bar{e} , f( \bar{e} )) $ and $ \alpha_0 = \alpha(x_0,f^*) $  \\
(1) & Iterate: & Compute $ x_{k+1}  := x_k + \smfrac{\alpha_k - 1 }{ \alpha_k \| P(g') \|^2} P(g') $, \,   where $ g' \in G( x_k', t_k') $. \\
&& Let $ (x_{k+1}', t_{k+1}') := \pi' (x_{k+1}, f^*) $ and $ \alpha_{k+1} :=  \alpha(x_{k+1}, f^*) $.  \vspace{1mm} 
\end{tabular} 

\noindent 
\hrulefill

\vspace{2mm}

\begin{thm}  \label{thm.ab} 
For the feasible sequence $ \{ x_k' \} $ generated by Algorithm B, and for $ 0 < \epsilon < 1 $,
 \begin{align}
 &   \ell  \geq 4 \, \left( \frac{\bar{D}}{\hat{r}} \right)^2 \, \cdot \,      \left( \, \frac{4}{ 3  } \left( \frac{1 - \epsilon }{ \epsilon } \right)^2 + 4 \left(   \frac{1 - \epsilon}{\epsilon} \right)  \right.  \nonumber \\
& \qquad \qquad \qquad  \qquad  \qquad  \qquad   \left.  + \log_2 \left( \frac{1 - \epsilon }{ \epsilon } \right) + \log_2 \left( \frac{\bar{D}}{\hat{r}} \right)   \, + 1 \, \right)   \nonumber    \\ &  \qquad        \quad \Rightarrow \quad \min_{k \leq \ell } \,   \frac{f(x_k') - f^*}{\hat{f} - f^*} \, \leq \, \epsilon \; , \label{eqn.af} \end{align} 
 \end{thm} 
\vspace{2mm}

When $ S $ is polyhedral and $ f $ is piecewise linear, then in addition to the implication (\ref{eqn.af}), for some constants $ C_1 $ and $ C_2 $ there holds
\begin{equation}  \label{eqn.ag} 
     \ell \geq C_1 \log(1/\epsilon ) + C_2 \quad \Rightarrow \quad \min_{k \leq \ell} \frac{f(x_k') - f^*}{\hat{f} - f^*} \leq \epsilon \; , \end{equation}
     i.e., linear convergence.  The constants, however, are highly dependent on $ \mathrm{Feas} $ and  $ f $, making the linear convergence mostly a curiosity. \vspace{1mm}

In the following four sections, focus is exclusively on convex conic linear optimization problems. In the first of these sections (\S \ref{sect.b}),  the equivalent problem is explained and the key theory is developed.  Perhaps most surprising is that the transformation to an equivalent problem is simple and the theory is elementary. \vspace{1mm}

Representative algorithms for the conic setting are developed and analyzed in Sections~\ref{sect.d} and \ref{sect.e}, first in the ideal case where the optimal value is known, and then generally. (Groundwork for the development and analysis is laid in \S \ref{sect.c}) \vspace{1mm}

In Section~\ref{sect.f}, the general convex optimization problem (\ref{eqn.aa}) is recast into conic form, to which algorithms and analyses from earlier sections can be applied.  The consequent complexity results are stated directly in terms of the original optimization problem (\ref{eqn.aa}), but the algorithms remain partially abstract in that some key computations are not expressed directly in terms of the original problem.  \vspace{1mm}

Finally, in Section \ref{sect.g}, the remaining ties to the original optimization problem are established, allowing the algorithms applied to the conic recasting to instead be entirely expressed in terms of the original problem, resulting in Algorithms A and B above.  The paper closes with the proofs of Theorems~\ref{thm.aa}  and \ref{thm.ab}, and the proof of the claim regarding linear convergence (i.e., (\ref{eqn.ag})).  
\vspace{1mm}

 We do not discuss how to compute appropriate input $ \bar{e} $ when such a point exists but is unknown,  because at present we do not see with high generality how to {\em  cleanly} apply the framework to accomplish this, let alone see how to apply the framework to compute a ``good'' choice for $ \bar{e} $ from among the points in $ \int(S \cap \dom(f)) $.  (Still, compared to the traditional literature where it is assumed $ \mathrm{Feas} \subset \int(\dom(f)) $ and (arbitrary) points can readily be orthogonally projected onto $ \mathrm{Feas} $, assuming $ \bar{e} $ is known and line-searches can be done seems to us considerably less restrictive.)        \vspace{1mm}

 The core of this paper, \S \S \ref{sect.b}--\ref{sect.f},    is virtually identical with our arXiv posting \cite{renegar2015framework}, a paper focused on convex conic linear optimization problems rather than general convex optimization problems, and having only at the end a brief discussion of the relevance to general convex optimization.  (Motivated by reviewers' comments, we decided in revising the work for publication to give primary focus to general convex optimization, while still utilizing the conic setting as the natural venue for developing the equivalent problem to which traditional subgradient methods are applied.)   Closely related to \cite{renegar2015framework} is recent  work of Freund and Lu \cite{freund2015new}, who develop first-order methods for convex optimization problems in which the objective function has a known lower bound and is assumed to satisfy a Lipschitz condition. Their perspective provides an interesting juxtaposition to the conic-oriented presentation in \cite{renegar2015framework}.  Also related is \cite{renegar2015accelerated}, in which accelerated methods for hyperbolic programming are presented.
 \vspace{1mm}

We close the introduction by emphasizing the intent is to develop a framework, not to prescribe specific algorithms.  Different subgradient methods applied in the framework yield different algorithms.  Changing how a general optimization problem is recast into conic form also results in a different algorithm.  
\vspace{2mm}

\section{{\bf Key Theory}}  \label{sect.b}

The key theory is elementary and yet has been overlooked in the literature, a blind spot.  \vspace{1mm}

Continue to let $ {\mathcal E} $ denote a finite-dimensional Euclidean space.  \vspace{1mm}

Let $ {\mathcal K}  \subset {\mathcal E} $ be a proper, closed, convex cone with nonempty interior. \vspace{1mm}

Fix a vector $ e \in \int(  {\mathcal K} ) $.  We refer to $ e $ as the ``distinguished direction.''   For each $ x \in {\mathcal E} $, let 
\[   \lambda_{\min}(x) :=  \inf \{ \lambda \mid x - \lambda \, e \notin {\mathcal K}  \} \; , \]
that is, the scalar $ \lambda $ for which $ x - \lambda e $ lies in the boundary of $ {\mathcal K}  $.  (Existence and uniqueness of $ \lambda_{\min}(x) $ follows from $ e \in \int({\mathcal K} ) \neq {\mathcal E} $ and the assumption that $ {\mathcal K}  $ is a closed, convex cone.)  \vspace{1mm}

If, for example, $ {\mathcal E} = \Sym $ ($ n \times n $ symmetric matrices), $ {\mathcal K}  = \Symp $ (cone of positive semidefinite matrices), and $ e = I $ (the identity), then $ \lambda_{\min}(X) $ is the minimum eigenvalue of $ X $. \vspace{1mm}

On the other hand, if $ {\mathcal K}  = \realsnp $ (non-negative orthant) and  $ e $ is a vector with all positive coordinates, then $ \lambda_{\min}(x) = \min_j x_j/e_j $ for $ x \in \mathbb{R}^n $. Clearly, the value of $ \lambda_{\min}(x) $ depends on the distinguished direction $ e $ (a fact the reader should keep in mind since the notation does not reflect the dependence).  \vspace{1mm}

Obviously, $ {\mathcal K}  = \{ x \mid \lambda_{\min}(x) \geq 0 \} $ and $ \int({\mathcal K} ) = \{ x \mid \lambda_{\min}(x) > 0 \} $. Also,
\begin{equation}  \label{eqn.ba}
  \lambda_{\min}(sx + te) = s \, \lambda_{\min}(x) + t \quad \textrm{for all $ x \in {\mathcal E}$  and scalars $ s \geq 0 $,  $ t $} \; .  
  \end{equation}

Let 
\[ \bar{{\mathcal B}} := \{ v \in {\mathcal E} \mid e + v \in {\mathcal K} \textrm{ and }   e - v \in {\mathcal K}  \} \; , \]
a closed, centrally-symmetric, convex set with nonempty interior. 
Define the gauge (\cite[\S 15]{rockafellar1970convex})  on $ {\mathcal E} $ by 
\[ \| u \|_{\infty} := \inf  \{ t \geq 0 \mid u = tv \textrm{ for some $ v \in \bar{{\mathcal B}} $} \} \; , \]
which is easily shown to be a norm if and only if $ {\mathcal K} $ is a pointed cone. 
Let $ \bar{B}_{\infty}(x,r) $ denote the closed ``ball'' centered at $ x $ and of radius $ r $. Clearly, $ \bar{B}_{\infty}(0,1) = \bar{{\mathcal B}} $, and $ \bar{B}_{\infty}(e,1) $ is the largest subset of $ {\mathcal K}  $ that has symmetry point $ e $, i.e., for each $ v $, either both points $ e + v $ and $ e - v$ are in the set, or neither point is in the set.      
\vspace{1mm}

\begin{prop}  \label{prop.ba}
  The function $ x \mapsto \lambda_{\min}(x) $ is concave and Lipschitz continuous: 
\[    | \lambda_{\min}(x) - \lambda_{\min}(y) | \leq \| x - y \|_{\infty} \quad \textrm{for all $ x,y \in {\mathcal E} $} \; .   \]  
\end{prop} 
\noindent {\bf Proof:} Concavity follows easily from the convexity of $ {\mathcal K}  $, so we focus on establishing Lipschitz continuity.  \vspace{1mm}

Let $ x, y \in {\mathcal E} $.  According to (\ref{eqn.ba}), the difference $ \lambda_{\min}(x + te) - \lambda_{\min}(y+te) $ is independent of $ t $, and of course so is the quantity $ \| (x + te) - (y + te) \|_{\infty} \; . $ Consequently, in proving the Lipschitz continuity, we may assume $ x $ lies in the boundary of $ {\mathcal K}  $, that is, we may assume $ \lambda_{\min}(x) = 0 $. The goal, then, is to prove
\begin{equation}  \label{eqn.bb}
    | \lambda_{\min}(x + v) | \leq \| v \|_{ \infty} \quad \textrm{for all $ v \in {\mathcal E} $} \; . \end{equation} 
  
We consider two cases.  First assume $ x + v $ does not lie in the interior of $ {\mathcal K}  $, that is, assume $ \lambda_{\min}(x+v) \leq 0 $. Then, to establish (\ref{eqn.bb}), it suffices to show $ \lambda_{\min}(x + v) \geq - \| v \|_{\infty} \; , $ that is, to show
\begin{equation}  \label{eqn.bc}
    x + v + \| v \|_{\infty} \, e \in {\mathcal K}  \; . 
    \end{equation} 
However, 
\begin{equation}  \label{eqn.bd}
   v + \| v \|_{\infty} \, e \in \bar{B}_{\infty}(\| v \|_{\infty} \, e, \| v \|_{\infty}) \subseteq {\mathcal K}  \; , 
   \end{equation} 
the set containment due to $ {\mathcal K}  $ being a cone and, by construction, $ \bar{B}_{\infty}(e,1) \subseteq {\mathcal K}   $. Since $ x \in {\mathcal K}  $ (indeed, $ x $ is in the boundary of $ {\mathcal K}  $), (\ref{eqn.bc})  follows. \vspace{1mm}

Now consider the case $ x + v \in {\mathcal K}  $, i.e., $ \lambda_{\min}(x+v) \geq 0 $. To establish (\ref{eqn.bb}), it suffices to show $ \lambda_{\min}(x + v) \leq \|v\|_{\infty} \; ,  $ that is, to show
\[                x + v - \|v\|_{\infty} \, e \notin \int({\mathcal K} ) \; . \]
Assume otherwise, that is, assume 
\[    x = w + \| v \|_{\infty} \, e - v \quad \textrm{for some $ w \in \int({\mathcal K} ) $} \; . \]
Since $ \| v \|_{\infty} \, e - v \in {\mathcal K}  $ (by the set containment on the right of (\ref{eqn.bd})), it then follows that $ x \in \int({\mathcal K} ) $, a contradiction to $ x $ lying in the boundary of $ {\mathcal K}  $.    \hfill $ \Box $
 \vspace{3mm}
 
 In this section, the inner product is used primarily for expressing a conic optimization problem.  To allow a distinction between inner products, we denote evaluation for the ``modeling inner product'' on a pair $ u,v \in {\mathcal E} $ by $ u \cdot v $, whereas in later sections we denote -- as was done implicitly in the introduction -- evaluation for the ``computational inner product'' by $ \lin u, v \rin $. \vspace{1mm}

Let  $ \mathrm{Affine} \subseteq {\mathcal E} $ be an affine space, i.e., the translate of a subspace. For fixed $ c \in {\mathcal E} $, consider the conic program
\[  
  \left. \begin{array}{rl}
\inf & c \cdot x  \\
\textrm{s.t.} & x \in \mathrm{Affine}  \\
  & x \in {\mathcal K}  \; .  \end{array} \right\} \cp  \]  
Assume $ \mathrm{Affine} \cap \int({\mathcal K} ) $ -- the set of strictly feasible points -- is nonempty. Let $ z^* $ denote the optimal value. \vspace{1mm}

Assume $ c $ is not orthogonal to the subspace of which $ \mathrm{Affine} $ is a translate, since otherwise all feasible points are optimal.  This assumption implies that all optimal solutions for CP lie in the boundary of $ {\mathcal K}  $. \vspace{1mm}

 Fix a strictly feasible point, $ e $.  The point $ e $ serves as the distinguished direction. \vspace{1mm}

For scalars $ z \in \mathbb{R}  $, let
\begin{gather*}  \mathrm{Affine}_{z} := \{ x \in \mathrm{Affine} \mid   c \cdot x = z \}, \textrm{ and}   \\
\textrm{let $ {\mathcal L} $ denote the subspace of which these affine spaces are translates.} 
\end{gather*}

Presently we show that for any choice of $ z $   satisfying $ z < c \cdot e \; ,  $ CP can be easily transformed into an equivalent optimization problem in which the only constraint is $ x \in \mathrm{Affine}_{z} \; . $ We make a simple observation.   \vspace{1mm}

\begin{lemma}  \label{lem.bb} 
Assume $ \mathrm{CP} $ has bounded optimal value.  

$ \textrm{~} $ \qquad \qquad  \qquad  If $ x \in \mathrm{Affine} $ satisfies $ c \cdot x < c \cdot e $, then $ \lambdamin(x) < 1 \; . $ 
\end{lemma}
\noindent {\bf Proof:}  It follows from (\ref{eqn.ba}) that if  $ \lambda_{\min}(x) \geq 1 $, then $ e + t( x - e) $ is feasible for all $ t \geq 0 $.  As the function $ t \mapsto c \cdot  \big( e + t(x- e) \big) $ is strictly decreasing (because $ c \cdot x < c \cdot e $), this implies CP has unbounded optimal value, contrary to assumption.~$ \Box $
 \vspace{3mm}
 
For  $ x \in {\mathcal E}  $ satisfying $ \lambda_{\min}(x) < 1 $, let $ \pi(x) $ denote the point where the half-line beginning at $ e $ in direction $ x - e $ intersects the boundary of $ {\mathcal K}  $:
\begin{equation}  \label{eqn.bda}
   \pi(x) :=  e + \smfrac{1}{1 - \lambda_{\min}(x)} (x - e)    
   \end{equation} 
(to verify correctness of the expression, observe (\ref{eqn.ba}) implies $ \lambda_{\min}(\pi(x)) = 0 $). 
We refer to $ \pi(x) $ as the ``radial projection'' of $ x $. \vspace{1mm}

 The proof of the following result is straightforward, but because it is central to our development, we label the result as a theorem. \vspace{1mm}

\begin{thm} \label{thm.bc}
Let $ z $ be any value satisfying \, $   z < c \cdot e \; . $  If $ x^* $ solves
\begin{equation}  \label{eqn.be}
  \begin{array}{rl}
   \sup & \lambdamin(x) \\
 \mathrm{s.t.} &  x \in \mathrm{Affine}_{z} \; ,  \end{array} 
    \end{equation} 
then $ \pi( x^* ) $ is optimal for $ \mathrm{CP} $. Conversely, if $ \pi^*  $ is optimal for $ \mathrm{CP} $, then $ x^* :=   e + \frac{c \cdot e - z }{c \cdot e - z^* } ( \pi^* - e) $ is optimal for (\ref{eqn.be}), and $ \pi^*  = \pi( x^* ) $.  
\end{thm}
\noindent {\bf Proof:} Fix a value $ z $ satisfying $ z < c \cdot e  $. It is easily proven from the convexity of $ {\mathcal K}  $ that $ x \mapsto \pi(x) $ gives a one-to-one map from $ \mathrm{Affine}_{z} $      onto 
\begin{equation}  \label{eqn.bf}
  \{ \pi \in \mathrm{Affine} \cap  \mathrm{bdy}({\mathcal K} )  \mid  c \cdot \pi < c \cdot e \} \; , 
  \end{equation}
where $ \mathrm{bdy}({\mathcal K} )  $ denotes the boundary of $ {\mathcal K}  $. \vspace{1mm}

For $ x \in \mathrm{Affine}_{z} \; , $ the CP objective value of $ \pi(x) $ is 
\begin{align}
  c \cdot \pi(x) & = c \cdot \big( e+ \smfrac{1}{1 - \lambda_{ \min}(x)} (x - e) \big) \nonumber \\
                    & = c \cdot e + \smfrac{1}{1 - \lambda_{ \min }(x)} ( z - c \cdot e ) \; , \label{eqn.bg}
                    \end{align}
 a strictly-decreasing function of $ \lambda_{\min}(x) $.                    
 Since the map $ x \mapsto \pi(x) $ is a bijection between $ \mathrm{Affine}_{z} $ and the  set (\ref{eqn.bf}), the theorem readily follows. \hfill $ \Box $
 \vspace{3mm}
 
CP has been transformed into an equivalent linearly-constrained maximization problem with concave, Lipschitz-continuous objective function. Virtually any subgradient method --  rather, {\em  supgradient} method -- can be applied to this problem, the main cost per iteration being in computing a supgradient and projecting it onto the subspace $ {\mathcal L}  $.   \vspace{1mm}

For illustration, we digress to interpret the implications of the development thus far for the linear program
\[ 
  \left. \begin{array}{rl}
 \min_{x \in \mathbb{R}^n} & c^T x \\
\textrm{s.t.} & Ax = b \\
 & x \geq 0 \; ,  \end{array}   \right\} \, \mathrm{LP}  
\] 
assuming $ e = {\bf 1} $ (the vector of all ones), in which case $ \lambda_{\min}(x) = \min_j x_j \; , $ and $ \| \, \, \|_{\infty} $ is the $ \ell_{\infty} $ norm, i.e., $ \| v \|_{\infty} = \max_j |v_j| $.  Let the number of rows of $ A $ be $ m \geq 1 $, assume the rows are linearly independent, and assume $ c $ is not a linear combination of the rows (otherwise all feasible points are optimal). \vspace{1mm}

For any scalar $ z < c^T {\bf 1} $, Theorem~\ref{thm.bc} asserts that LP is equivalent to   
\begin{equation} \label{eqn.bh}
  \begin{array}{rl}
    \max_x & \min_j x_j \\
     \textrm{s.t.} & Ax = b \\
           & c^T x = z \; , \end{array} \end{equation} 
in that when $ x $ is feasible for (\ref{eqn.bh}), $ x $ is optimal if and only if the projection \\ $ \pi(x) = \mathbf{1}  + \smfrac{1}{1 - \min_j x_j} (x - \mathbf{1} ) $ is optimal for LP. 
  The setup is shown schematically in the following figure:

$ \textrm{~} $ \quad \qquad \quad \qquad  \qquad   \includegraphics[scale=.26]{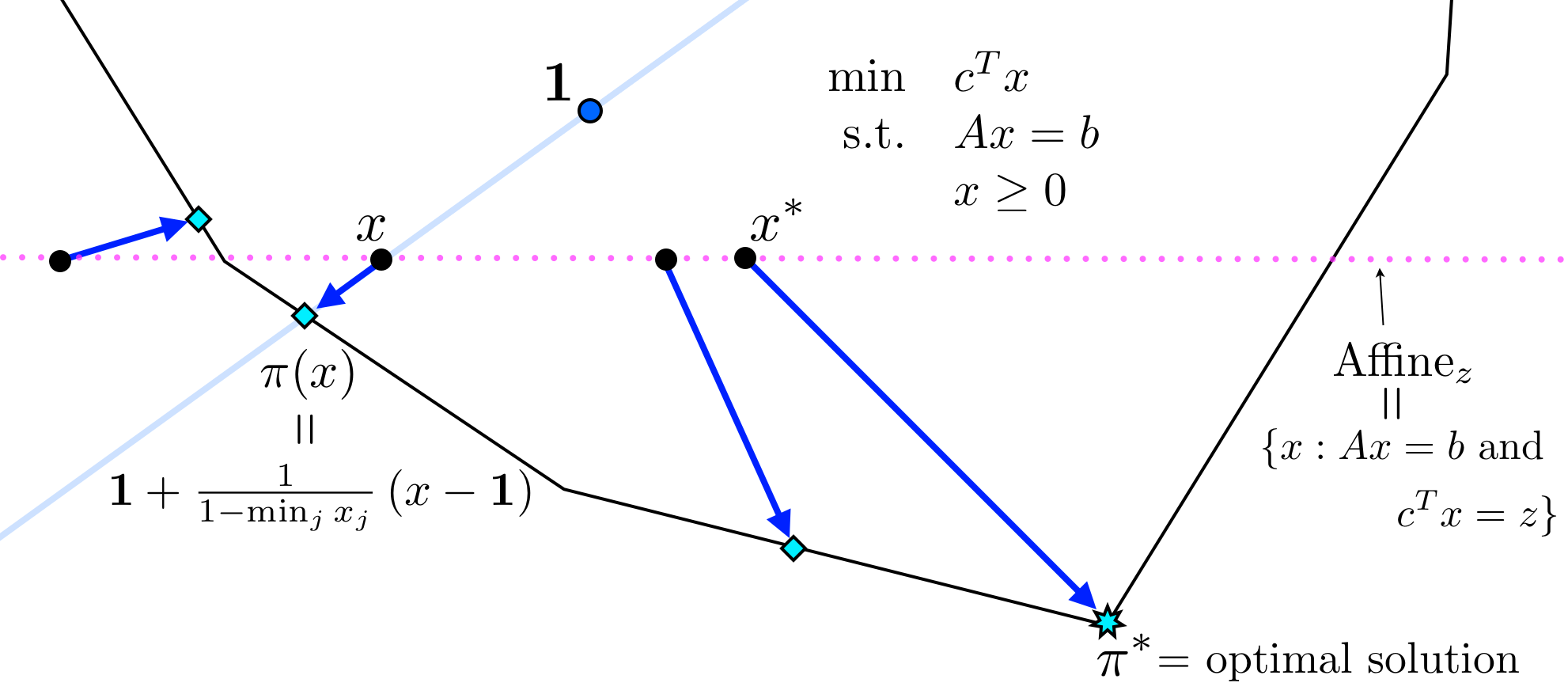}
\vspace{4mm}

Proposition~\ref{prop.ba}  asserts that, as is obviously true, $ x \mapsto \min_j x_j $ is $ \ell_{\infty} $-Lipschitz continuous with constant 1. Consequently, the function also is $ \ell_2$-Lipschitz continuous with constant 1, as is relevant if supgradient methods rely on the standard inner product in computing supgradients.    \vspace{1mm}

With respect to the standard inner product, the supgradients of $ x \mapsto \min_j x_j $ at $ x $ are  the convex combinations of the standard basis vectors $ e(k) $ for which $ x_k = \min_j x_j $. Consequently, the projected supgradients at $ x $ are the convex combinations of the vectors $ \bar{P}_k $ for which $ x_k = \min_j x_j $, where $ \bar{P}_k $ is the $ k^{th} $ column of the matrix projecting $ \mathbb{R}^n $ onto  the nullspace of $ \bar{A} = \left[ \begin{smallmatrix}  A \\ c^T \end{smallmatrix} \right] $, that is  
\[ 
   \bar{P}  := I- \bar{A}^T (\bar{A} \,   \bar{A}^T)^{-1} \bar{A} \; .   \] 

If $ m \ll n $, then $ \bar{P} $ is not computed in its entirety, but instead the matrix $ \bar{M}  = ( \bar{A} \bar{A}^T)^{-1} $ is formed as a preprocessing step, at cost $ O(m^2 n) $ (the inverse exists because the set consisting of $ c $ and the rows of $ A $ has been assumed linearly independent). Then, for any iterate $ x $ and an index $ k $ satisfying $ x_k = \min_j x_j $, the projected  supgradient $ \bar{P}_k $ is computed according to
 \[  u = \bar{M} \, \bar{A}_k \quad \rightarrow \quad  v = \bar{A}^T u \quad \rightarrow \quad \bar{P}_k = e(k) - v \; , \]
for a cost of $ \, O(m^2 \, + \, \# \mathrm{non\_zero\_entries\_in\_} A   )  $ per iteration.  \vspace{1mm}

Before returning to the general theory, we note that if the choices are $ {\mathcal E} = \Sym $, $ {\mathcal K}  = \Symp $ and $ e = I $ (and thus $ \lambda_{\min}(X) $ is the minimum eigenvalue of $ X $), then with respect to the trace inner product, the supgradients at $ X $ for the function $ X \mapsto \lambda_{\min}(X) $ are the convex combinations of the matrices $ v v^T $, where $ Xv = \lambda_{\min}(X) v $ and $ \| v \|_2 = 1 $.   \vspace{1mm}

{\em  Assume, henceforth, that CP has at least one optimal solution, and that $ z $ is a fixed scalar satisfying \, $ z < c \cdot e $.}   Then the equivalent problem (\ref{eqn.be}) has at least one optimal solution.     Let $ x^*_{z} $ denote any of the optimal solutions for the equivalent problem, and recall
$ z^* $ denotes the optimal value of CP. A useful characterization of the optimal value for the equivalent problem is easily provided. \vspace{1mm}

\begin{lemma}  \label{lem.bd}
 \[   \lambda_{\min}(x^*_{z }) = \frac{z - z^* }{c \cdot e - 
 z^* } \] 
  \end{lemma}
\noindent {\bf Proof:}   
By Theorem~\ref{thm.bc}, $ \pi(x^*_{z}) $ is optimal for CP -- in particular, $ c \cdot  \pi(x^*_{z})  = z^* $.   Thus, according to (\ref{eqn.bg}),
\[  z^* = c \cdot e + \smfrac{1}{1 - \lambda_{\min}(x^*_{z })} \, ( z - c \cdot e) \; . \]
Rearrangement completes the proof.      \hfill $ \Box $
 \vspace{3mm}
 
We focus on the goal of computing a point $ \pi $ which is feasible for CP and has better objective value than  $ e$ in that 
\begin{equation}  \label{eqn.bi}
     \frac{c \cdot \pi - z^* }{c \cdot e - z^* } \leq \epsilon \; \; , 
     \end{equation} 
where $ 0 < \epsilon < 1 $.  Thus, for the problem of primary interest, CP, the focus is on relative improvement in the objective value. 
\vspace{1mm}

 The following proposition provides a useful characterization of the accuracy needed in approximately solving the CP-equivalent problem (\ref{eqn.be}) so as to ensure that for the computed point $ x $, the projection $ \pi = \pi(x) $ satisfies (\ref{eqn.bi}). \vspace{1mm}

\begin{prop}  \label{prop.be}
If $ x \in \mathrm{Affine}_z $ and $ 0 < \epsilon < 1 $,  then 
\begin{align*} 
   & \frac{ c \cdot  \pi(x)  - z^* }{c \cdot e - z^* } \, \leq \,  \epsilon   
    \\ & \qquad   \qquad  \qquad   \textrm{if and only if}    \\ & \qquad  \qquad  \qquad  \qquad   
    \lambdamin( x^*_{ z}  ) - \lambdamin( x) \, \leq \, \frac{\epsilon }{1 - \epsilon } \, \, \frac{c \cdot e - z }{\, \, \, c \cdot e - z^* } \; .  
\end{align*} 
\end{prop}
\noindent {\bf Proof:}  Assume $ x \in \mathrm{Affine}_{z} $.  For $ y = x $ and $ y = x^*_{z} \; , $  we have the equality (\ref{eqn.bg}), that is,  
\[      c \cdot  \pi(y)  = c \cdot e + \smfrac{1}{1 - \lambda_{\min}(y)} ( z - c \cdot e ) \; . 
\]
Thus,
\begin{align*}
    \frac{ c \cdot  \pi(x) - z^* }{c \cdot e - z^* } & =   \frac{ c \cdot  \pi(x)  - c \cdot  \pi(x^*_{z} )  }{c \cdot e - c \cdot  \pi(x^*_{z} )  } \\
    & = \frac{ \smfrac{1}{1 - \lambdamin(x)} -  \smfrac{1}{1 - \lambdamin(x^*_{z} )} }{ - \smfrac{1}{1 - \lambdamin(x^*_{z} )}} \\  
    & = \frac{ \lambdamin(x^*_{z} ) - \lambdamin(x) }{ 1 - \lambdamin(x)} \; . 
    \end{align*}
Hence,
\begin{gather*}     
 \frac{c \cdot  \pi(x)  - z^* }{c \cdot e - z^* } \leq \epsilon \\
\Leftrightarrow \\ 
   \lambdamin(x^*_{z} ) - \lambdamin(x) \leq \epsilon \, ( 1 - \lambdamin(x)) \\
\Leftrightarrow \\
(1 - \epsilon ) ( \lambdamin(x^*_{z} ) - \lambdamin(x) ) \leq \epsilon ( 1 - \lambdamin(x^*_{z} )) \\
\Leftrightarrow \\ \lambdamin(x^*_{z} ) - \lambdamin(x) \leq \smfrac{\epsilon }{1 - \epsilon } ( 1 - \lambdamin(x^*_{z} )) \; .  
   \end{gather*}
Using Lemma~\ref{lem.bd}    to substitute for the rightmost occurrence of $ \lambda_{\min}(x^*_{z} ) $ completes the proof. \hfill $ \Box $
 \vspace{3mm}

\section{ {\bf Groundwork for Algorithm Design and Analysis}} \label{sect.c}

In Sections \ref{sect.d} and \ref{sect.e}, we show how the key theory leads to algorithms and complexity results regarding the solution of the conic program CP. In this section, groundwork is laid for the development in those and subsequent sections.
\vspace{1mm}

  Continue to assume CP has an optimal solution, denote the optimal value by $ z^* $, and let $ e$ be a strictly feasible point, the distinguished direction.  
Given $ \epsilon > 0 $ and a value $ z $ satisfying $ z < c \cdot e $, the approach is to apply supgradient methods to approximately solve
\begin{equation}  \label{eqn.ca}
 \begin{array}{rl}
 \max & \lambda_{ \min}(x) \\
\textrm{s.t.} & x \in \mathrm{Affine}_{z}   \; ,  \end{array} 
      \end{equation}
where by ``approximately solve'' we mean that  $ x \in \mathrm{Affine}_{z}  $ is computed for which 
\[ 
\lambda_{\min}(x^*_{z}) - \lambda_{ \min}(x) \,  \leq \,  \frac{\epsilon}{1 - \epsilon } \, \, \frac{c \cdot e - z   }{\, \, \, c \cdot e - z^* }  \; . 
\] 
Indeed, according to Proposition~\ref{prop.be}, the projection $ \pi = \pi(x) $ will then satisfy
\[ 
\frac{c \cdot \pi - z^* }{c \cdot e - z^* } \leq \epsilon  \; . 
\] 
\vspace{1mm}

Not until results regarding CP are applied to general convex optimization -- Sections \ref{sect.f} and \ref{sect.g}   -- do we require a characterization of supgradients of the function $ x \mapsto \lambdamin(x) $, but as the characterization is interesting in itself, we present it as the first piece of groundwork laid in this section. \vspace{1mm}

Let $ P_{{\mathcal L}} $ denote orthogonal projection onto $ {\mathcal L} $ with respect to $ \langle \; , \; \rangle $ (the computational inner product), and let $ \hat{\partial} \lambda_{\min}(x) $ denote the supdifferential at $ x $. The following proposition relates supdifferentials to normal cones (normal with respect to $ \langle \; , \; \rangle $). \vspace{1mm}

\begin{prop}  \label{prop.ca} 
For all $ x \in {\mathcal E} $, 
  \[      \hat{\partial}\lambda_{\min}(x) = \{ v  \mid  -v \in N_{{\mathcal K}}( x - \lambda_{\min}(x) e) \textrm{ and } \lin e, v \rin = 1 \} \; . \]
\end{prop}
\noindent {\bf Proof:}  We know 
\begin{equation}  \label{eqn.cb} 
  \lambda_{\min}(y + re) = \lambda_{\min}(y) + t \quad \textrm{for all $ y \in {\mathcal E} $, $ t \in \mathbb{R} $}  \; , 
  \end{equation}
  from which follows for every $ x $ and $ t $,
\begin{equation}   \label{eqn.cc} 
  \hat{\partial}\lambda_{\min}(x + te)  = \hat{\partial} \lambda_{\min}(x) 
  \end{equation} 
and
\begin{equation}  \label{eqn.cd} 
          v \in \hat{\partial} \lambda_{\min}(x) \quad \Rightarrow \quad \lin e, v \rin = 1 \; .
\end{equation}
Due to (\ref{eqn.cc}), in proving the lemma we may assume $ \lambdamin(x) = 0 $, i.e., $ x \in \mathrm{bdy}({\mathcal K}) $.  
\vspace{1mm}

Since $ \lambdamin(x) = 0 $, and the value of $ \lambdamin $ is non-negative on $ {\mathcal K} $, 
\[  g \in \hat{\partial} \lambdamin(x) \quad \Rightarrow \quad \forall y \in {\mathcal K}, \,  \lin g, y - x \rin \geq 0 \; . \] 
Thus, with (\ref{eqn.cd}), we have 
\begin{equation}  \label{eqn.ce} 
   \hat{\partial}\lambdamin(x) \subseteq \{ v \mid \lin e, v \rin = 1 \textrm{ and }    - v \in N_{{\mathcal K}}(x) \} \; . 
   \end{equation}

On the other hand, if $ -v \in N_{{\mathcal K}}(x) $, then 
\[  y \in \mathrm{bdy}({\mathcal K}) \quad \Rightarrow \quad      0 = \lambdamin(y) \leq \lambdamin(x) + \lin v, y - x \rin \]
(using $ \lambdamin(x) = 0 $).   Thus, if in addition, $ \lin e, v \rin = 1 $, 
\[  t \in \mathbb{R} \textrm{ and }  y \in \mathrm{bdy}({\mathcal K}) \quad \Rightarrow \quad \lambdamin(y + te) \leq \lambdamin(x) + \lin v, (y + te) - x \rin \]   
(using \ref{eqn.cb}). Since $ {\mathcal E} = \{ y + te \mid y \in \mathrm{bdy}( {\mathcal K}) \textrm{ and } t \in \mathbb{R} \} $, the reverse inclusion to (\ref{eqn.ce})   thus holds, establishing the proposition. \hfill $ \Box $
 \vspace{3mm}

Let $ \| \, \, \| $ be the norm associated with $ \langle \; , \; \rangle $, and for $ z \in \mathbb{R} $, let
\[  M_z := \sup \left\{ \smfrac{| \lambda_{\min}(x) - \lambda_{\min}(y)|}{ \| x - y \| } \mid x,y \in \mathrm{Affine}_z  \textrm{ and } x \neq y \right\} \; , \]
the Lipschitz constant for the map $ x \mapsto \lambda_{\min}(x) $ restricted to $ \mathrm{Affine}_z $.  Proposition~\ref{prop.ba}  implies $ M_z $ is well-defined (finite), although unlike the Lipschitz constant for the norm appearing there (i.e., $ \| \, \, \|_{\infty} $), $ M_z $ might exceed 1, depending on $ \| \, \, \| $. The next piece of groundwork to be laid is a geometric characterization of an upper bound on $ M_z $. \vspace{1mm}

We claim the values $ M_z $ are identical for all $ z $.  To see why, consider that for $ z_1, z_2 < c \cdot e \; , $  a bijection from $ \mathrm{Affine}_{z_1} $ onto $ \mathrm{Affine}_{z_2} $ is provided by the map
\[  x \mapsto y(x) := \smfrac{z_2 - z_1}{c \cdot e -  z_1} e + \smfrac{c \cdot e - z_2}{c \cdot e - z_1}  x \; . \]
Observe, using (\ref{eqn.ba}),  
\[  \lambda_{\min}(y(x)) = \smfrac{c \cdot e - z_2}{c \cdot e - z_1} \lambda_{\min}(x) + \smfrac{z_2 - z_1}{c \cdot e - z_1} \; , \]
and thus 
\[  \lambda_{\min}(y(x)) - \lambda_{\min}(y( \bar{x})) = \smfrac{c \cdot e - z_2}{c \cdot e - z_1} \left(   \lambda_{\min}(x) - \lambda_{\min}( \bar{x}) \right) \quad \textrm{for $ x, \bar{x} \in \mathrm{Affine}_{z_1} $} \; . \]
Since, additionally, $ \| y(x) - y( \bar{x}) \| = \smfrac{c \cdot e - z_2}{c \cdot e - z_1} \| x - \bar{x} \| $, it is immediate that the values $ M_z $ are identical for all $ z < c \cdot e $.  A simple continuity argument then implies this value is equal to $ M_{c \cdot e} $. Analogous reasoning shows $ M_z = M_{ c \cdot e} $ for all $ z > c \cdot e $.  In all, $ M_z $ is independent of $ z $, as claimed.  \vspace{1mm}

Let $ M $ denote the common value, i.e., $ M = M_z $ for all $ z $. \vspace{1mm}

The following proposition can be useful in modeling and in choosing the computational inner product. This result plays a central role in our complexity bounds for general convex optimization. \vspace{1mm}

Let $ \bar{B}(e,r) := \{ x \mid \| x - e \| \leq r \} $.  

\begin{prop} \label{prop.cb} $ M \leq 1/ r_e \; , $ where  $ r_e :=  \max \{ r \mid \bar{B}(e,r) \cap \mathrm{Affine}_{c \cdot e} \subseteq {\mathcal K}  \} $ 
\end{prop}
\noindent {\bf Proof:} According to Proposition~\ref{prop.ba}, 
\[   | \lambda_{\min}(x) - \lambda_{\min}(y)| \leq \| x - y \|_{\infty} \quad \textrm{for all $ x, y $} \; . \]
Consequently, it suffices to show $ \| x - y \| \geq r_e \| x - y \|_{\infty} $ for all $ x,y \in \mathrm{Affine}_{c \cdot e} \; , $ i.e., it suffices to show for all $  v \in {\mathcal L} $ that $ \| v \| \geq r_e \| v \|_{\infty} $. \vspace{1mm}

 However, according to the discussion just prior to Proposition~\ref{prop.ba},  $ \bar{B}_{\infty}(e,1) $ is the largest set which both is contained in $ {\mathcal K} $ and has symmetry point $ e $, from which follows that $ \bar{B}_{\infty}(e,1) \cap \mathrm{Affine}_{c \cdot e} $ is the largest set which is both contained in $ {\mathcal K}  \cap \mathrm{Affine}_{c \cdot e} $ and has symmetry point $ e $.  Hence 
 \[  \bar{B}(e, r_e) \cap \mathrm{Affine}_{ c \cdot e} \subseteq \bar{B}_{\infty}(e,1) \cap \mathrm{Affine}_{c \cdot e} \; , \]
implying $ \| v \| \geq r_e \| v \|_{\infty} $ for all $ v \in {\mathcal L} $. \hfill $ \Box $ \vspace{3mm}

Towards considering specific supgradient methods, we recall the following standard and elementary result, rephrased for our setting:

\begin{lemma}   \label{lem.cc} 
  Assume $ z \in \mathbb{R} $, $ x, y \in \mathrm{Affine}_z $ and $ g \in \hat{\partial}\lambdamin(x) $. 

For all scalars $ \alpha  $,
\[ 
   \| ( x + \alpha P_{{\mathcal L}} (g)) - y \|^2 \leq \| x - y \|^2 - 2 \alpha  \left(   \lambda_{\min}(y) - \lambda_{\min}(x) \right) + \alpha^2 \| P_{{\mathcal L}} (g) \|^2 \; . \] 
\end{lemma} 
\noindent {\bf Proof:}  Letting $ \tilde{g} := P_{{\mathcal L}} (g )$,   simply observe
\begin{align*}
\| ( x + \alpha \tilde{g} ) - y \|^2 & = \| x - y \|^2 + 2 \alpha \lin \tilde{g} , x - y \rin + \alpha^2 \| \tilde{g}  \|^2 \\
& =  \| x - y \|^2 - 2 \alpha \lin g, y - x \rin + \alpha^2 \| \tilde{g}  \|^2 \quad \textrm{(by $ x - y \in {\mathcal L} $)} \\
& \leq \| x - y \|^2 - 2 \alpha \left(   \lambda_{\min}(y) - \lambda_{\min}(x) \right)  + \alpha^2 \| \tilde{g}  \|^2 \; , 
\end{align*}
the inequality due to concavity of the map $ x \mapsto \lambda_{\min}(x) \; . $ \hfill $ \Box $ \vspace{2mm} 

\section{{\bf Algorithm 1: When the Optimal Value is Known}} \label{sect.d} 

Knowing $ z^* $ is not an entirely implausible situation.  For example, if strict feasibility holds for a primal conic program and for its dual
\[ \begin{array}{rl}
          \min & \bar{c}^T x \\
          \textrm{s.t.} & \bar{A}x = \bar{b} \\
     & x \in \bar{{\mathcal K} } \end{array} \qquad  \begin{array}{rl}
          \max & \bar{b}^T y \\
           \textrm{s.t.} & \bar{A}^T y + s = \bar{c} \\
                    & s \in \bar{{\mathcal K} }^* \; , \end{array} \]
then the combined primal-dual conic program is known to have optimal value equal to zero:
\[   \begin{array}{rl}
       \min & \bar{c}^T x - \bar{b}^T y \\
       \textrm{s.t.} & \bar{A} x = \bar{b} \\
                    & \bar{A}^T y + s = \bar{c} \\
                    & (x,s) \in \bar{{\mathcal K} } \times \bar{{\mathcal K} }^* \; . \end{array} \]
        
 \noindent 
\hrulefill
                    
\noindent {\bf Algorithm 1} \vspace{1mm} 

\noindent 
\noindent 
\begin{tabular}{lll} 
(0) &Input: & $ z^* $, the optimal value of CP, \\
&&  $ e $, a strictly feasible point for CP, and \\
&&    $ \bar{x}  \in \mathrm{Affine} $ satisfying $ c \cdot \bar{x}  < c \cdot e \; . $  \\
& Initialize: & Let $ x_0 = e + \frac{c \cdot e - z^*}{c \cdot e - c \cdot \bar{x} }(\bar{x} -e) $ \,  \,  (thus, $ c \cdot x_0 = z^* $), \\
&&      and let $ \pi_0 = \pi(x_0) $ \, \,  ($ = \pi(\bar{x} ) $). \\
(1) & Iterate: &  Compute $ x_{k+1} = x_k - \smfrac{\lambda_{\min}(x_k)}{\| P_{{\mathcal L}} ( g_k ) \|^2} P_{{\mathcal L}} ( g_k ) $, \,  where $ g_k \in \hat{\partial}\lambdamin(x_k) $.   \\
&&  Let $ \pi_{k+1} = \pi(x_{k+1}) $. \vspace{1mm} 
\end{tabular} 

\noindent 
\hrulefill
\vspace{2mm}

All of the iterates $ x_k $ lie in $ \mathrm{Affine}_{z^*} $, and hence, $ \lambda_{\min}(x_k) \leq 0 $, with equality if and only if $ x_k $ is feasible (and optimal) for CP. \vspace{1mm}

For all scalars $ z < c \cdot e $ and for $ x \in \mathrm{Affine}_z \; ,  $ define
\[  
   \mathrm{dist}_z(x) :=  \min \{ \| x - x_z^* \|  \mid x_z^* \textrm{ is optimal for  (\ref{eqn.ca})} \} \; . \] 

\begin{prop}  \label{prop.da} 
The iterates for Algorithm 1 satisfy
\[ 
 \max \{ \lambda_{\min}(x_k)  \mid k = \ell, \ldots, \ell + m \} \geq -M \, \mathrm{dist}_{z^*}(x_{\ell})/ \sqrt{m+1} \; . 
\] 
\end{prop} 
\noindent {\bf Proof:}  Letting $ \tilde{g}_k := P_{{\mathcal L}} (g_k) $,  Lemma~\ref{lem.cc}  implies
\begin{align*} 
   & \mathrm{dist}_{z^*}(x_{k+1})^2 \\ & \leq  \mathrm{dist}_{z^*}(x_k)^2 - 2 (-\lambda_{\min}(x_k)/\| \tilde{g}_k \|^2) \,  ( 0 - \lambda_{\min}(x_k)) \, + \, (\lambda_{\min}(x_k)/\| \tilde{g}_k \|)^2 \\  
    & = \mathrm{dist}_{z^*}(x_k)^2 - (\lambda_{\min}(x_k)/\| \tilde{g}_k \|)^2 \; , 
    \end{align*}
and thus by induction (and using $ \| \tilde{g}_k \| \leq M $),
\begin{align*}
\mathrm{dist}_{z^*}(x_{\ell + m + 1 })^2 & \leq \mathrm{dist}_{z^*}(x_{\ell})^2 - \sum_{k= \ell }^{\ell + m} ( \lambda_{\min}(x_k)/M)^2  \\
& \leq  \mathrm{dist}_{z^*}(x_{\ell})^2 - \smfrac{m + 1}{M^2}  \min \{ \lambda_{\min}(x_k)^2  \mid k = \ell, \ldots, \ell + m \} \; , 
\end{align*}
implying the proposition (keeping in mind $ \lambda_{\min}(x_k) \leq 0 $). \hfill $ \Box $
 \vspace{3mm}
  
We briefly digress to consider the case of $ {\mathcal K}  $ being polyhedral, where already an interesting result is easily proven. The following corollary is offered only as a curiosity, as the constants  typically are so large as to render the bound on $ \ell$ meaningless except for minuscule $ \epsilon$.
 \vspace{1mm}

\begin{cor}  \label{cor.db} 
Assume $ {\mathcal K}  $ is polyhedral. There exist constants $ C_1 $ and $ C_2 $ (dependent on CP, $ e $, $ \bar{x} $ and the computational inner product), such that for all $ 0 < \epsilon < 1 $, 
\[  \ell  \geq C_1 + C_2 \log(1/ \epsilon ) \quad \Rightarrow \quad \min_{k \leq \ell} \frac{c \cdot \pi_k - z^*}{c \cdot e - z^*} \,  \leq \, \epsilon \; . \]
\end{cor}
\vspace{2mm}

  For first-order methods, such a logarithmic bound in $ \epsilon $  was initially established by Gilpin, Pe\~{n}a and Sandholm \cite{gilpin2012first}. They did not assume an initial feasible point $ e $ was known, but neither did they require the computed approximate solution to be feasible (instead, constraint residuals were required to be small).  They relied on an accelerated gradient method, along with the smoothing technique of Nesterov \cite{nesterov2005smooth}.  As is the case for the above result, they assumed the optimal value of CP to be known apriori, and they restricted $ {\mathcal K}  $ to be polyhedral. \vspace{2mm}
  
The proof of the corollary depends on the following simple lemma. \vspace{1mm}        

\begin{lemma}  \label{lem.dc} 
For Algorithm 1, the iterates satisfy
\[  \frac{c \cdot \pi_k - z^*}{c \cdot e - z^*} = \frac{- \lambda_{\min}(x_k) }{1 - \lambda_{\min}(x_k)} \; . \]
\end{lemma}
\noindent {\bf Proof:}  Immediate from $ \pi_k = e + \frac{1}{1 - \lambda_{\min}(x_k)} (x_k - e) $ and $ c \cdot x_k = z^* $. \hfill $ \Box $
 \vspace{3mm}  
 
\noindent {\bf Proof of Corollary~\ref{cor.db}:} 
 With $ {\mathcal K}  $ being polyhedral, the concave function $ x \mapsto \lambda_{\min}(x) $ is piecewise linear, and thus there exists a positive constant $ C $ such that
\[      \mathrm{dist}_{z^*}(x) \leq - C \, \lambda_{\min}(x) \quad \textrm{for all $ x \in \mathrm{Affine}_{z^*} $} \; . \]
Then Proposition~\ref{prop.da}  gives
\[  
 \max  \{ \lambda_{\min}(x_k)  \mid k = \ell, \ldots, \ell + m \} \geq C \, M \, \lambda_{\min}(x_{\ell})/ \sqrt{m+1} \; , \]
from which follows 
 \[  \max  \{ \lambda_{\min}(x_k)  \mid k = \ell, \ldots, \ell + \lceil (2CM)^2 \rceil \} \geq \smfrac{1}{2} \lambda_{\min}(x_{\ell}) \; , \]
i.e., $ \lambda_{\min}(x_{\ell}) $ is ``halved'' within $ \lceil (2CM)^2 \rceil $ iterations.  The proof is easily completed using Lemma~\ref{lem.dc}.   \hfill $ \Box $ \vspace{3mm}

We now return to considering general convex cones $ {\mathcal K}  $. \vspace{1mm}

The iteration bound provided by Proposition~\ref{prop.da}  bears little obvious connection to the geometry of the conic program CP, except in that the constant $ M $ is related to the geometry by Proposition~\ref{prop.cb}.  The other constant -- $ \mathrm{dist}_{z^*}(x_k) $ -- does not at present have such a clear geometrical connection to CP.  We next observe a meaningful connection.  \vspace{1mm}

 The {\em level sets}  for CP are the sets
\[  \mathrm{Level}_{z} =   \mathrm{Affine}_{z} \cap {\mathcal K}   \; ,   \]
that is, the largest feasible sets for CP on which the objective function is constant\footnote{There is possibility of confusion here, as in the optimization literature, the terminology ``level set'' is often used for the portion of the feasible region on which the (convex) objective function does not exceed a specified value rather than -- as for us -- exactly equals the value.  Our terminology is consistent with the general mathematical literature, where the region on which a function does not exceed a specified value is referred to as a sublevel set, not a level set.}.  If $ z < z^* $, then $ \mathrm{Level}_{z} = \emptyset \; . $    \vspace{1mm}

If some level set is unbounded, then either CP has unbounded optimal value or can be made to have unbounded value with an arbitrarily small perturbation of $ c $.  Thus, in developing numerical optimization methods, it is natural to focus on the case that level sets for CP are bounded.  \vspace{1mm}  
     
For scalars $ z  $, define the diameter of $ \mathrm{Level}_z $ by  
\[    \diamval  :=  \sup  \{  \| x - y \| \mid x,y \in \mathrm{Level}_{z}   \} \; ,    \]
the diameter of $ \mathrm{Level}_{z} $.  If $ \mathrm{Level}_z = \emptyset $, let $ \diamval := - \infty  $. \vspace{1mm}

\begin{lemma}  \label{lem.dd} 
Assume $ x \in \mathrm{Affine}_{z^*} $, and let $ \pi = \pi(x) $. Then
       \[  \mathrm{dist}_{z^*}(x) \,  = \,  (1 - \lambda_{\min}(x)) \, \mathrm{dist}_{c \cdot \pi}(\pi) \,  = \, \frac{ \mathrm{dist}_{c \cdot \pi}(\pi)}{1 - \frac{c \cdot \pi - z^*}{c \cdot e - z^*}} \,  \leq \, \frac{ \mathrm{diam}_{c \cdot \pi}}{1 - \frac{c \cdot \pi - z^*}{c \cdot e - z^*}} \; .  \] 
 \end{lemma}
\noindent {\bf Proof:}  Since 
\begin{equation}  \label{eqn.da} 
  \pi = e + \smfrac{1}{1 - \lambda_{\min}(x)}(x-e) \; , \end{equation}
  Theorem~\ref{thm.bc}  implies that the maximizers of the map $ y \mapsto \lambda_{\min}(y) $ over $ \mathrm{Affine}_{ c \cdot \pi} $ are precisely the points of the form
\[   x_{c \cdot \pi}^* =   e + \smfrac{1}{1 - \lambda_{\min}(x)}(x_{z^*}^*-e) \; , \]
where $ x_{z^*}^* $ is a maximizer of the map when restricted to $ \mathrm{Affine}_{z^*} $ (i.e., is an optimal solution of CP).  Observing  
\[    \pi - x_{c \cdot \pi}^* = \smfrac{1}{1 - \lambda_{\min}(x)}  (x - x_{z^*}^*) \; , \]
it follows that 
\[   \mathrm{dist}_{c \cdot \pi}( \pi) =  \smfrac{1}{1 - \lambda_{\min}(x)}  \, \mathrm{dist}_{z^*}(x) \; , \]
establishing the first equality in the statement of the lemma.  The second equality then follows easily from (\ref{eqn.da}) and $ c \cdot x = z^* $. The inequality is due simply to $ \pi, x_{c \cdot \pi}^* \in \mathrm{Level}_{ c \cdot \pi} $, for all optimal solutions $ x_{c \cdot \pi}^* $ of the CP-equivalent problem (\ref{eqn.ca}) (with $ z = c \cdot \pi $).   \hfill $ \Box $
 \vspace{3mm}
 
For scalars $ z $, define
\[  \mathrm{Diam}_z := \max \{ \mathrm{diam}_{z'} \mid z' \leq z \} \; , \]
the ``horizontal diameter'' of the sublevel set consisting of points $ x $ that are feasible for CP and satisfy $ c \cdot x \leq z $.  For $ z^* < z < c \cdot e $, the  value $ \mathrm{Diam}_{z} $ can be thought of as a kind of condition number for CP, because $  \mathrm{Diam}_{z} $ being large is an indication that the optimal value for CP is relatively sensitive to perturbations in the objective vector $ c $. \vspace{1mm}

For $ z^* \leq  z < c \cdot e $, define
\[     \mathrm{Dist}_z := \sup \{ \mathrm{dist}_{z'}(x) \mid z' \leq z \textrm{ and } x \in \mathrm{Level}_{z'} \}  \; . \]
Clearly, there holds the relation
\[    \mathrm{Dist}_z \leq \mathrm{Diam}_z \; , \]
and hence if the ``condition number'' $ \mathrm{Diam}_z $ is only of modest size, so is the value $ \mathrm{Dist}_z $.  \vspace{1mm}

Following is our main result for Algorithm 1. By substituting $  \mathrm{Diam}_{c \cdot \pi_0} $ for $  \mathrm{Dist}_{ c \cdot \pi_0} $, and $ 1/r_e$ for $ M $ (where $ r_e $ is as in Proposition~\ref{prop.cb}), the statement of the theorem becomes phrased in terms clearly reflecting the geometry of CP. \vspace{1mm}

\begin{thm}  \label{thm.de} 
Assume $ 0 < \epsilon < \frac{c \cdot \pi_0 - z^*}{c \cdot e - z^*} $, where $ \pi_0 = \pi(x_0) $ is the initial $ \mathrm{CP} $-feasible point for Algorithm 1 (i.e., assume $ \pi_0 $ does not itself satisfy the desired accuracy).  Then 
\begin{align*}
 &   \ell  \geq (2M \, \mathrm{Dist}_{c \cdot \pi_0})^2 \, \, \left( \, \frac{4}{ 3  } \left( \frac{1 - \epsilon }{ \epsilon } \right)^2 + 4 \left(   \frac{1 - \epsilon}{\epsilon} \right)  \right. \\
& \qquad \qquad \qquad  \qquad  \qquad  \qquad   \left. + \log_2 \left( \frac{\frac{c \cdot \pi_0 - z^*}{c \cdot e - z^*} }{\epsilon} \right) + \log_2 \left( \frac{1 - \epsilon }{1 - \frac{c \cdot \pi_0 - z^*}{c \cdot e - z^*} } \right) \, + 1 \, \right)       \\ &  \qquad        \quad \Rightarrow \quad \min_{k \leq \ell } \,   \frac{c \cdot \pi_k - z^*}{c \cdot e - z^*} \, \leq \, \epsilon \; . \end{align*} 
\end{thm}
\noindent {\bf Proof:}   To ease notation, let $  \lambda_k := \lambda_{\min}(x_k) \; . $ \vspace{1mm}

 Let $ k_0 = 0 $ and recursively define $ k_{i+1} $ to be the first index for which $ \lambda_{k_{i+1}} \geq \lambda_{k_i}/2 $ (keeping in mind $ \lambda_k \leq 0 $ for all $ k $).  Proposition~\ref{prop.da} implies  
\begin{align}
   k_{i+1} - k_i + 1   & \leq \left( \frac{2M \, \mathrm{dist}_{z^*}(x_{k_i})}{\lambda_{k_i}} \right)^2 \nonumber \\
   & = \left( \, 2 M \, \mathrm{dist}_{c \cdot \pi_{k_i}}(\pi_{k_i}) \, \frac{ 1 - \lambda_{k_i}}{ \lambda_{k_i}} \,  \right)^2 \quad \textrm{(by Lemma~\ref{lem.dd})} \nonumber \\ 
  & \leq   \left( \, 2 M \, \mathrm{Dist}_{c \cdot \pi_0} \, \frac{ 1 - \lambda_{k_i}}{ \lambda_{k_i}} \,  \right)^2 \; ,  \label{eqn.db} 
  \end{align}
  where the final inequality is due to $ c \cdot \pi_{k_i} $ ($ i = 0, 1, \ldots $) being a decreasing sequence (using Lemma~\ref{lem.dc}). \vspace{1mm}
  
Let $ i' $ be the first sub-index for which $ \lambda_{k_{i'}} \geq - \epsilon/(1 - \epsilon) $. Lemma~\ref{lem.dc}  implies
\[   \frac{c \cdot \pi_{k_{i'}} - z^*}{ c \cdot e - z^* } \, \leq  \, \epsilon \; . \]
Thus, to prove the theorem, it suffices to show $ \ell  = k_{i'} $ satisfies the inequality in the statement of the theorem.  \vspace{1mm}

Note $ i' > 0 $ (because, by assumption, $ \epsilon < \frac{c \cdot \pi_0 - z^*}{c \cdot e - z^*} $). Observe, then, 
\begin{align}
    i' &  < 1 +   \log_2 \left( \frac{ \lambda_0}{- \epsilon/(1 - \epsilon )} \right)  \nonumber \\
       &  = 1 +   \log_2 \left( \frac{\frac{c \cdot \pi_0 - z^*}{c \cdot e - z^*} }{ \epsilon} \right) + \log_2 \left( \frac{1 - \epsilon }{1 - \frac{c \cdot \pi_0 - z^*}{c \cdot e - z^*} } \right)    \label{eqn.dc} 
       \end{align}  
       (again using Lemma~\ref{lem.dc}).     
  \vspace{1mm}

Additionally,
\begin{align*}
k_{i'} & = \sum_{i=0}^{i'-1} k_{i+1} - k_i  \\ &
    \leq (2M \, \mathrm{Dist}_{c \cdot \pi_0})^2 \, \sum_{i=0}^{i'-1} \left( \frac{1 - \lambda_{k_i}}{\lambda_{k_i}} \right)^2 \quad \textrm{(by (\ref{eqn.db}))} \\
       & \leq (2M \, \mathrm{Dist}_{c \cdot \pi_0})^2 \,  \sum_{i=0}^{i'-1} \left( \frac{1 - 2^i \lambda_{k_{i'-1}}}{2^i \lambda_{k_{i'-1}}} \right)^2  \\
       & \leq  (2M \, \mathrm{Dist}_{c \cdot \pi_0})^2 \, \sum_{i=0}^{i'-1} \left( \frac{1 + 2^i \epsilon /(1 - \epsilon)}{2^i \epsilon /(1 - \epsilon) } \right)^2 \\
     &  = (2M \, \mathrm{Dist}_{c \cdot \pi_0})^2 \, \sum_{i=0}^{i'-1} \left( 1 + \frac{1}{2^i} \, \frac{1 - \epsilon}{\epsilon } \right)^2 \\    
& \leq (2M \, \mathrm{Dist}_{c \cdot \pi_0})^2 \, \left(  i' + 4  \frac{1- \epsilon }{ \epsilon }  + \frac{4}{3} \left( \frac{1- \epsilon }{ \epsilon } \right)^2  \right) \; .    
\end{align*}
Using (\ref{eqn.dc})  to substitute for $ i' $ completes the proof. \hfill $ \Box $

\section{ {\bf Algorithm 2: When The Optimal Value Is Unknown}} \label{sect.e}

For the second algorithm, we discard the requirement of knowing $ z^* $.  Now  $ \epsilon $ (the desired relative-accuracy) is required as input. \vspace{2mm} 

\noindent 
\hrulefill
 
\noindent {\bf Algorithm 2} \vspace{1mm}

\noindent   
\begin{tabular}{lll}
(0) & Input: & $ 0 < \epsilon < 1 $ ,\\
&&  $ e $, a strictly feasible point for CP, and \\
&&   $ \bar{x}  \in \mathrm{Affine} $ satisfying $ c \cdot \bar{x}  < c \cdot e $.  \\
& Initialize: & $ x_0 = \pi_0 = \pi(\bar{x} ) $ \\
(1) &Iterate: & Compute $ \tilde{x}_{k+1} :=  x_k + \smfrac{ \epsilon }{2 \| P_{{\mathcal L}} g_k \|^2} P_{{\mathcal L}} g_k $, \, where $ g_k \in \hat{\partial}\lambdamin(x_k) $.  \\
&& Let $ \pi_{k+1} :=  \pi(\tilde{x}_{k+1}) \; . $ \\
&& If $ c \cdot (e - \pi_{k+1}) \geq \smfrac{4}{3}  \, c \cdot (e - \tilde{x}_{k+1}) $, let $ x_{k+1} = \pi_{k+1} $; \\
&& \qquad     else, let $ x_{k+1} = \tilde{x}_{k+1} \; . $   \vspace{1mm}
\end{tabular} 

\noindent 
\hrulefill
\vspace{2mm}

Unsurprisingly, the iteration bound we obtain for Algorithm 2 is worse than the result for Algorithm 1, but perhaps surprisingly, the bound is not excessively worse, in that the factor for $ 1/ \epsilon^2 $ is essentially unchanged (it's the factor for $ 1/\epsilon $ that increases, although typically not by a large amount). \vspace{2mm}

\begin{thm}  \label{thm.ea} 
 Assume $ 0 < \epsilon < \frac{c \cdot \pi_0 - z^*}{c \cdot e - z^*}  $.  For the iterates of Algorithm 2,
\begin{align*}
 &   \ell  \geq 8 \, (M \, \mathrm{Dist}_{c \cdot \pi_0} )^2 \, \left( \, \frac{1}{\epsilon^2} \, + \, \frac{1}{\epsilon} \, \log_{4/3} \left( \frac{ 1  }{1 - \frac{c \cdot \pi_0 - z^*}{c \cdot e - z^*}} \right)  \,  \right)              \\ &  \qquad        \quad \Rightarrow \quad \min_{k \leq \ell} \,   \frac{c \cdot \pi_k - z^*}{c \cdot e - z^*} \, \leq \, \epsilon \; . \end{align*} 
\end{thm}
\noindent {\bf Proof:} In order to distinguish the iterates obtained by projecting to the boundary, we record a notationally-embellished rendition of the algorithm which introduces a distinction between ``inner iterations'' and ``outer iterations'':
\vspace{2mm}

\noindent 
\hrulefill
 
\noindent  {\em Algorithm 2  (notationally-embellished version):} \vspace{1mm}

\noindent  
\begin{tabular}{l}
\begin{tabular}{lll}
(0) & Input: & $ 0 < \epsilon < 1 \; , $ $ e $ and $ \bar{x} $.  \\ 
& Initialize: & $ y_{1,0} = \pi(\bar{x} ) \; , $  \\
&&  $ i = 1$ (outer iteration counter), \\
&&    $ j = 0 $ (inner iteration counter). 
\end{tabular} \\
\begin{tabular}{ll}
(1) & Compute $ y_{i,j+1} = y_{i,j} + \smfrac{\epsilon }{2 \| \tilde{g}_{i,j} \|^2} \, \tilde{g}_{i,j} \; , $  \\
& \qquad  \qquad  \qquad  \qquad  \qquad  \qquad  \quad    where $ \tilde{g}_{i,j} = P_{{\mathcal L}}g_{i,j} $ and $ g_{i,j} \in \hat{\partial}\lambdamin(y_{i,j})   $. \\
(2) & If $ c \cdot (e - \pi(y_{i,j+1})) \geq \smfrac{4}{3} \, c \cdot (e - y_{i,j+1} ) \; , $ \\
  & \qquad  then let \, $ y_{i+1,0} = \pi(y_{i,j+1}) $, \,  $ i \leftarrow i + 1  $ \, and  \, $  j \leftarrow 0 \; $;\\
&  \quad   else, let \, $ j \leftarrow j + 1 \; . $  \\
(3) & Go to step 1.  \vspace{1mm}
\end{tabular} 
\end{tabular}  

\noindent 
\hrulefill
\vspace{2mm}

For each outer iteration $ i $, all of the iterates $ y_{i,j} $ have the same objective value.  Denote the value by $ z_i $.  Obviously, $ z_1 $ is equal to the value $ c \cdot \pi_0 $ appearing in the statement of the theorem.  Let
\[ \mathrm{Dist} := \mathrm{Dist}_{c \cdot \pi_0} \, = \mathrm{Dist}_{z_1} \; . \]

Step 2 ensures     
\begin{equation}  \label{eqn.ea} 
    c \cdot e - z_{i+1} \geq \smfrac{4}{3}  ( c \cdot e - z_i)  \; . 
\end{equation}  
Thus, $ z_1, z_2, \ldots $ is a strictly decreasing sequence.  Consequently, as $ y_{i,0} \in \mathrm{Level}_{z_i} $, we have $  \mathrm{dist}_{z_i}(y_{i,0}) \leq \mathrm{Dist} $ for all $ i $.  \vspace{1mm}

From (\ref{eqn.ea})  we find for scalars $ \delta > 0 $ that
\[  
 \frac{c \cdot e - z_{i+1}}{ c \cdot e - z^*} <  \delta \quad \Rightarrow \quad i < \log_{4/3} \left( \frac{\delta}{ \frac{c \cdot e - z_1}{c \cdot e - z^*}} \right) \,  = \log_{4/3} \left( \frac{\delta}{ 1 - \frac{z_1 - z^*}{c \cdot e - z^*}} \right) \; , 
\]
and thus, for $ \epsilon < 1 $, 
\begin{equation}   \label{eqn.eb} 
  \frac{z_i - z^*}{c \cdot e - z^*} >  \epsilon \quad \Rightarrow \quad  i <  1 + \log_{4/3} \left( \frac{1 - \epsilon}{ 1 - \frac{z_1 - z^*}{c \cdot e - z^*}} \right) \; . 
  \end{equation}
Hence, if an outer iteration $ i $ fails to satisfy the inequality on the right, the initial inner iterate $ y_{i,0} $ fulfills the goal of finding a CP-feasible point $ \pi $  satisfying $ \frac{c \cdot \pi - z^*}{c \cdot e - z^*} \leq \epsilon $ (i.e., the algorithm has been successful no later than the start of outer iteration $ i $).  Also observe that (\ref{eqn.eb})  provides (letting $ \epsilon \downarrow 0 $) an upper bound on $ I $, the total number of outer iterations:
\begin{equation}  \label{eqn.ec} 
          I \leq 1 + \log_{4/3} \left( \frac{1}{ 1 - \frac{z_1 - z^*}{c \cdot e - z^*}} \right) \; . 
  \end{equation}
  \vspace{1mm}

For $ i = 1, \ldots, I $,  let $ J_i $ denote the number of inner iterates computed during outer iteration $ i $, that is, $ J_i $ is the largest value $ j $ for which $ y_{i,j} $ is computed. Clearly, $ J_I = \infty $, whereas $ J_1, \ldots, J_{I-1} $ are finite. \vspace{1mm}

To ease notation, let $ \lambda_{i,j} := \lambda_{\min}(y_{i,j}) $, and let $ \lambda_i^* := \lambda_{\min}(x_{z_i}^*) $, the optimal value of
\[ \begin{array}{rl}
    \max & \lambda_{\min}(x) \\
     \textrm{s.t.} & x \in \mathrm{Affine}_{z_i} \; . \end{array} \]
According to Lemma~\ref{lem.bd},
\begin{equation}  \label{eqn.ed} 
\lambda_i^* = \frac{z_i - z^*}{c \cdot e - z^*} \; . 
\end{equation}
It is thus valid, for example, to substitue $ \lambda_i^* $ for $ \frac{z_i - z^*}{c \cdot e - z^*} $ in (\ref{eqn.eb}).   Additionally, (\ref{eqn.ed})  implies (\ref{eqn.ea})  to be equivalent to
\begin{equation}  \label{eqn.ee} 
   1 - \lambda_{i+1}^* \geq \smfrac{4}{3} ( 1 - \lambda_i^* ) \; . 
   \end{equation} 
         
For any point $ y $, we have $ \pi(y) = e + \smfrac{1}{1 - \lambda_{\min}(y)} ( y - e) $,  and thus,
\[   \frac{ c \cdot e - c \cdot \pi(y)}{c \cdot e - c \cdot y} = \frac{1}{1 - \lambda_{\min}(y)} \; . \]
Hence, 
\[ \frac{c \cdot e - c \cdot \pi(y)}{c \cdot e - c \cdot y} \geq \frac{4}{3} \quad \Leftrightarrow \quad \lambda_{\min}(y) \geq 1/4 \; . \]
Consequently,   
\begin{equation}  \label{eqn.ef} 
          \lambda_{i,j} < 1/4 \, \, \textrm{ for $ j < J_i $} \; . 
         \end{equation} 
   
We use the following relation implied by Lemma~\ref{lem.cc}: 
\begin{equation}  \label{eqn.eg} 
 \mathrm{dist}_{z_i}(y_{i,j+1})^2 \leq  \mathrm{dist}_{z_i}(y_{i,j})^2 - \smfrac{\epsilon }{\| \tilde{g}_{i,j} \|^2 } ( \lambda_i^* - \lambda_{i,j} ) + (  \smfrac{\epsilon}{2 \| \tilde{g}_{i,j} \| })^2 \; .
 \end{equation}

We begin bounding the number of inner iterations by showing
 \begin{equation}  \label{eqn.eh} 
       \lambda_i^* \geq \max \{ \smfrac{1}{2} , \epsilon \}  \quad \Rightarrow \quad  J_i \leq \frac{8 (M \, \mathrm{Dist})^2}{ \epsilon } \; .   
       \end{equation}
Indeed, for $ j < J_i \; , $  
\begin{align*}
 &    - \epsilon \,  ( \lambda_i^* - \lambda_{i,j} ) + \smfrac{1}{4} \,  \epsilon^2 \\
& <  - \epsilon \, ( \max \{ \smfrac{1}{2}, \epsilon \}  - \smfrac{1}{4} ) + \smfrac{1}{4} \epsilon^2 \quad \textrm{(using (\ref{eqn.ef}))} \\
& = \min  \left\{ \smfrac{1}{4} ( \epsilon^2 - \epsilon), \smfrac{1}{4} \epsilon - \smfrac{3}{4}  \epsilon^2  \right\}   \\
& \leq \smfrac{3}{4} \,   \smfrac{1}{4} ( \epsilon^2 - \epsilon) + \smfrac{1}{4} ( \smfrac{1}{4} \epsilon - \smfrac{3}{4}  \epsilon^2 ) \\
& = - \smfrac{1}{8} \epsilon \; .     
\end{align*} 
Thus, according to (\ref{eqn.eg}), for $ j < J_i $,
\[  \mathrm{dist}_{z_i}(y_{i,j+1})^2 \leq  \mathrm{dist}_{z_i}(y_{i,j})^2 - \frac{ \epsilon}{8  M^2} \; , \]
inductively giving
\begin{align*}
    \mathrm{dist}_{z_i}(y_{i,j+1})^2 & \leq \mathrm{dist}_{z_i}(y_{i,0})^2 - \frac{(j+1) \, \epsilon}{8 \,  M^2 } \\
        & \leq \mathrm{Dist}^2 - \frac{ (j+1) \, \epsilon}{ 8  M^2} \; . 
        \end{align*}
  The implication (\ref{eqn.eh})  immediately follows. \vspace{1mm}      

The theorem is now readily established in the case $ \epsilon \geq 1/2 $. Indeed, because of the identity (\ref{eqn.ed}), the quantity on the right of (\ref{eqn.eb})  provides an upper bound on the number of outer iterations $ i $ for which $ \lambda_i^* > \epsilon $, whereas the quantity on the right of (\ref{eqn.eh})  gives, assuming $ \epsilon \geq 1/2 $, an upper bound on the number of inner iterations for each of these outer iterations.  However, for the first outer iteration satisfying $ \lambda_i^* \leq \epsilon $, the initial iterate $ y_{i,0} $ ($ = \pi(y_{i,0})) $ itself achieves the desired accuracy $ \frac{c \cdot \pi - z^*}{c \cdot e - z^*} \leq \epsilon $. Thus, the total number of inner iterations made before the algorithm is successful is at most the product of the two quantities which  is seen not to exceed the iteration bound in the statement of the theorem (using $ \log_{4/3}(1- \epsilon) \leq  \log_{4/3}(1/2) < -1 $). \vspace{1mm}

Before considering the remaining case, $  \epsilon < 1/2 $, we establish a relation applying for all $ \epsilon $.  For any outer iteration $ i $  for which $ \lambda_i^* < 3/4 $, and for any $ 0 < \epsilon < 1 $,  let
\[  \widehat{J}_i := \left\lceil  \frac{1}{\smfrac{3}{4}  - \lambda_i^*} \,  \left( \frac{M \, \mathrm{Dist}}{  \epsilon } \right)^2 - 1 \right\rceil  \; . \]
We claim that either
\begin{equation}  \label{eqn.ei} 
    J_i \leq   \widehat{J}_i \quad \textrm{or} \quad  \min \left\{ \frac{c \cdot \pi(y_{i, j }) - z^*}{ c \cdot e - z^* } \mid j = 0, \ldots, \widehat{J}_i \right\} \, \leq \, \epsilon    \; . 
\end{equation}
Consequently, if $ J_i > \widehat{J}_i $,  the algorithm will achieve the goal of computing a point $ y $ satisfying $ \frac{c \cdot \pi(y) - z^*}{c \cdot e - z^*} \leq \epsilon  $  within $ \widehat{J}_i $ inner iterations during outer iteration $ i $. 
\vspace{1mm}

To establish (\ref{eqn.ei}), assume $ \widehat{J}_i <  J_i $ and yet the inequality on the right of (\ref{eqn.ei})  does not hold.  (We obtain a contradiction.) For every $ j \leq \widehat{J}_i $, Proposition~\ref{prop.be}  then implies 
\begin{align*}
 \lambda_i^* -  \lambda_{i,j} & > \epsilon \, \,  \frac{c \cdot e - z_i}{c \cdot e - z^*}  \\
                 & =  (1 - \lambda_i^*) \, \epsilon  \quad \textrm{(by (\ref{eqn.ed}))} \; ,
\end{align*}                  
 and hence, using (\ref{eqn.eg}),
 \[  \mathrm{dist}_{z_i}(y_{i,j+1})^2 <  \mathrm{dist}_{z_i}(y_{i,j})^2   - ( \smfrac{3}{4}  - \lambda_i^*) \,  \left( \epsilon / M \right)^2 \; , \]
from which inductively follows
\begin{align*}
   \mathrm{dist}_{z_i}( y_{i,\widehat{J}_i+1})^2 & <  \mathrm{dist}_{z_i}( y_{i,0})^2 - (\widehat{J}_i +1) \, (\smfrac{3}{4} - \lambda_i^* ) \, ( \epsilon/M)^2 \\
& \leq \mathrm{Dist}^2  - (\widehat{J}_i +1) \, (\smfrac{3}{4} - \lambda_i^* ) \, ( \epsilon/M)^2 \\
& \leq    0 \; , 
\end{align*}
a contradiction. The claim is established. \vspace{1mm}

Assume $ \epsilon < 1/2 $, the case remaining to be considered. \vspace{1mm}

As each outer iteration $ i $ satisfying $ \lambda_i^* \geq 1/2 $ has only finitely many inner iterations, there must be at least one outer iteration $ i $ satisfying $ \lambda_i^* < 1/2 $.  
Let $ i $ be the first outer iteration for which $ \lambda_i^* < 1/2 $.  From (\ref{eqn.ec})  and (\ref{eqn.eh}), the total number of inner iterations made before reaching outer iteration $ i $ is at most
\begin{equation}  \label{eqn.ej} 
\frac{8 (M \, \mathrm{Dist})^2}{ \epsilon } \, \,  \log_{4/3} \left( \frac{1}{ 1 - \frac{z_1 - z^*}{c \cdot e - z^*}} \right)   \; .    
  \end{equation} 

According to (\ref{eqn.ei}), during outer iteration $ i $, the algorithm either achieves its goal within $ \widehat{J}_i $ inner iterations, or the algorithm makes no more than $ \widehat{J}_i $ inner iterations before starting a new outer iteration.  Assume the latter case. Then, for outer iteration $ i + 1 $, the algorithm either achieves its goal within $ \widehat{J}_{i+1} $ inner iterations, or the algorithm makes no more than $ \widehat{J}_{i+1} $ inner iterations before starting a new outer iteration. Assume the latter case. In iteration $ i + 2 $, the algorithm definitely achieves its goal within $ \widehat{J}_{i+2} $ inner iterations, because there cannot be a subsequent outer iteration due, by (\ref{eqn.ee}), to
\[   \smfrac{4}{3} (1 - \lambda_{i+2}^*) \geq \left( \smfrac{4}{3} \right)^3 (1 - \lambda_i^*) >  \left( \smfrac{4}{3} \right)^3 \smfrac{1}{2} > 1 \; . \]

The total number of inner iterations made before the algorithm achieves its goal is thus bounded by the sum of the quantity (\ref{eqn.ej})   and 
\begin{align*}
          &    \widehat{J}_i + \widehat{J}_{i+1} + \widehat{J}_{i+2} \\
&  <  \left(  \frac{1}{(1-\frac{1}{2}) - \frac{1}{4}  } + \frac{1}{\frac{4}{3} \, (1 - \smfrac{1}{2}) - \smfrac{1}{4}  } + \frac{1}{\frac{4}{3} \,  \frac{4}{3} \, (1 - \smfrac{1}{2}) - \smfrac{1}{4}}  \right) \, \, \left( \frac{M \, \mathrm{Dist}}{ \epsilon } \right)^2 \\
& \qquad   \textrm{(using $ \smfrac{3}{4} - \lambda_j^* = (1 - \lambda_j^*) - \smfrac{1}{4} $)} \\ 
& < 8 \,  \left( \frac{M \, \mathrm{Dist}}{ \epsilon } \right)^2 \; , 
\end{align*}
completing the proof of the theorem. \hfill $ \Box $

\section{{\bf  Application to General Convex Optimization}} \label{sect.f}

We now return to the setting of Section 1, considering optimization problems of the form
\begin{equation}   \label{eqn.fa}
  \begin{array}{rl}  \min & f(x) \\
                           \textrm{s.t.} & x \in \mathrm{Feas} \; , \end{array} 
                           \end{equation}  
where $ f: {\mathcal E} \rightarrow (- \infty , \infty ] $  is an extended-valued and lower-semicontinuous convex function, where $ \mathrm{Feas} = \{ x \in S \mid Ax = b \} $, with $ S $ being a closed convex set, and where there is known a point $ \bar{e} \in \mathrm{int}(S \cap \dom(f)) $.
\vspace{1mm}

In this section we recast (\ref{eqn.fa})  into conic form, then interpret the complexity results obtained by applying Algorithms 1 and 2 to the conic refomulation. Before proceeding, we recall notation from \S \ref{sect.a}. \vspace{1mm}

Let $ \langle \; , \; \rangle $ be the (computational) inner product on $ {\mathcal E} $, and $ \| \, \, \| $ the associated norm. Let $ P $ be the linear operator orthogonally projecting $ {\mathcal E} $ onto the kernel of $ A $. \vspace{1mm}

Recall  $ \bar{D}  $ denotes the diameter of the sublevel set $ \{ x \in \mathrm{Feas} \mid f(x) \leq f( \bar{e} ) \} $.  The diameter is assumed to be finite, implying the optimal value $ f^* $ is attained at some feasible point.   \vspace{1mm}

Recall $ \hat{f} $ is a user-chosen scalar satisfying $ \hat{f} > f( \bar{e} ) $ (hence, $ ( \bar{e}, \hat{f}) \in \int(\epi(f)) $), and recall
\begin{equation} \label{eqn.fb} 
  \hat{r}  := \sup \{  r \mid \| x - \bar{e}  \| \leq r \textrm{ and }   Ax = b \,  \, \Rightarrow  \, \, x \in S \textrm{ and }   f(x) \leq \hat{f} \}  \; , 
  \end{equation} 
a positive scalar (because $ \bar{e}  \in \mathrm{int}\big( S \cap  \mathrm{dom}(f)\big) $).   \vspace{1mm}

The values $ \bar{D} $ and $ \hat{r} $ are not assumed to be known, but do appear in complexity results. \vspace{1mm}

 For later reference, observe the convexity of $ f $ implies
\[     f( \bar{e} )  \leq \smfrac{\hat{r} }{ \bar{D}  + \hat{r} } f^* + \smfrac{\bar{D}  }{\bar{D}  + \hat{r} } \hat{f}   \; , \]
which in turn implies
\begin{equation}  \label{eqn.fc} 
\frac{1}{1 - \frac{f( \bar{e} ) - f^*}{ \hat{f} - f^*}} \, \leq \, 1 + \bar{D}/\hat{r}     \; . 
\end{equation}

As $ S $ and $ \epi(f) $ are closed and convex,  there exist closed, convex cones $ {\mathcal K}_1, {\mathcal K}_2 \subseteq {\mathcal E} \times \mathbb{R} \times \mathbb{R} $ for which
\[  S \times \mathbb{R}  = \{ (x,t)  \mid (x,1, t) \in {\mathcal K}_1 \}  \quad \textrm{and} \quad      \mathrm{epi}(f) = \{ (x,t) \mid (x,1,t) \in {\mathcal K}_2 \} \; . \]
Letting $ {\mathcal K} := {\mathcal K}_1 \cap {\mathcal K}_2 $,    
clearly the optimization problem (\ref{eqn.da})  is equivalent to
\begin{equation}  \label{eqn.fd}
  \begin{array}{rl}
 \min_{x,s,t} & t \\
      \textrm{s.t.} & Ax = b \\
                    & s = 1 \\
   & (x,s,t) \in {\mathcal K} \; , \end{array} \end{equation} 
and has the same optimal value, $ f^* $.  The conic program (\ref{eqn.fd})  is of the same form as CP, the focus of preceding sections. Clearly,
\[  \mathrm{Affine} = \{ (x,1,t) \mid Ax = b \} \; , \textrm{ and for scalars $ z $, } \mathrm{Affine}_z = \{ (x,1,z) \mid Ax = b \} \; . \]

For distinguished direction, choose $ e = ( \bar{e} , 1, \hat{f} \, ) $, which clearly lies in the interiors of $ {\mathcal K}_1 $ and $ {\mathcal K}_2 $, and thus lies in $ \int( {\mathcal K}) $.   This distinguished direction, along with the cone $ {\mathcal K} $, determines the map $ (x,s,t) \mapsto \lambda_{\min}(x,s,t) $ on $ {\mathcal E} \times \mathbb{R} \times \mathbb{R} $. \vspace{1mm}  

For $ z < f^* $, the conic problem \ref{eqn.fd}  -- and hence the problem \ref{eqn.fa}  -- is, by Theorem~\ref{thm.bc}, equivalent to 
\begin{equation}   \label{eqn.fe} 
 \begin{array}{rl}
         \max_{x,s,t} & \lambda_{\min}(x,s,t) \\
          \textrm{s.t.} & (x,s,t) \in \mathrm{Affine}_z \; ,
          \end{array} 
          \end{equation} 
It is to this problem that supgradient methods are applied.
\vspace{1mm}

Choose the computational inner product on $ {\mathcal E} \times \mathbb{R} \times \mathbb{R} $ to to be 
\[ \blin (x_1, s_1, t_1), (x_2, s_2, t_2) \brin := \lin x_1, x_2 \rin + s_1 s_2 + t_1 t_2 \; ,   \]
in which case
\[  P_{{\mathcal L}}(x,s,t) = (P(x),0,0) \; . \]  
  \vspace{1mm}

Observe for all scalars $ z $,
\begin{equation}  \label{eqn.ff}
 \mathrm{Level}_{z} = \{ (x,1, z ) \mid x \in \mathrm{Feas} \textrm{ and }   f(x) \leq z   \} \; . \end{equation}
 
Recall that the Lipschitz constant $ M $ for the map $ (x,s,t) \mapsto \lambda_{\min}(x,s,t) $ restricted to $ \mathrm{Affine}_z $ is independent of $ z $. Since by  
(\ref{eqn.fb}), 
\[  \mathrm{Affine}_{\hat{f}} \cap \bar{B}(e,\hat{r} ) \subseteq {\mathcal K} \; , \]
Proposition~\ref{prop.cb} implies
\begin{equation}  \label{eqn.fg} 
              M \leq 1/ \hat{r} \; . 
              \end{equation}

Choose the input $ \bar{x}  $ to Algorithms 1 and 2 as $ \bar{x}  = ( \bar{e} , 1, f( \bar{e} )) $, which clearly is feasible for the conic program (\ref{eqn.fd}).  Note that (\ref{eqn.ff})  then implies the horizontal diameter of the relevant sublevel set for the conic program satisfies 
\begin{equation}  \label{eqn.fh} 
   \mathrm{Diam}_{f(\bar{e} )} =  \bar{D} \; .    
   \end{equation}
 
Algorithm 2 requires input $ 0 < \epsilon < 1 $, but not $ f^* $. 
Applying Algorithm 2 results in a sequence of iterates $ (x_k, 1, t_k) $ for which the projections $ (x_k',1,t_k') := \pi(x_k,1,t_k) $ satisfy $ x_k' \in \mathrm{Feas} $ and $ f(x_k') \leq t_k' $ (simply because $ (x_k',1,t_k') $ is feasible for the conic program (\ref{eqn.fd})). \vspace{1mm}

Since $ (x_0,1,t_0) = (\bar{e} ,1,f(\bar{e} )) \in \bdy({\mathcal K}) $, we have $  \pi(x_0, 1, t_0) = (x_0,1,t_0) $ -- in particular, the objective value of $ \pi(x_0,1,t_0) $ is $ f(\bar{e})  $.  Consequently, the sequence of points $ x_k' $ not only lie in $ \mathrm{Feas} $, but by Theorem~\ref{thm.ea}  satisfies   
\begin{align}  
 & \ell  \geq 8 \, \left(  \frac{\bar{D}}{\hat{r}} \right)^2    \, \left( \, \frac{1}{\epsilon^2} \, + \, \frac{1}{\epsilon} \, \log_{4/3} \left( 1 + \frac{\bar{D}}{\hat{r}} \right)       \, \right)   \label{eqn.fi}   \\
&  \quad                \quad \Rightarrow \quad \min_{k \leq \ell } \,   \frac{f(x_k') - f^*}{\hat{f} - f^*} \, \leq \, \epsilon \; ,   \label{eqn.fj} 
 \end{align}
where for (\ref{eqn.fi})  we have used (\ref{eqn.fc}), (\ref{eqn.fg}), (\ref{eqn.fh}), and for (\ref{eqn.fj})  have used $ f(x_k') \leq t_k' $. \vspace{1mm}
               
 For Algorithm 1, which requires input $ f^* $ but not $ \epsilon $, points $ (x_k', 1, t_k') $ are generated for which the sequence $ \{ x_k' \} $ is feasible and, by Theorem~\ref{thm.de}, satisfies
 \begin{align}
 &  
  \ell  \geq 4 \, \left( \frac{\bar{D}}{\hat{r}} \right)^2 \, \cdot \,      \left( \, \frac{4}{ 3  } \left( \frac{1 - \epsilon }{ \epsilon } \right)^2 + 4 \left(   \frac{1 - \epsilon}{\epsilon} \right)  \right.  \label{eqn.fk} \\
& \qquad \qquad \qquad  \qquad  \qquad  \qquad   \left.  + \log_2 \left( \frac{1 - \epsilon }{ \epsilon } \right) + \log_2 \left( \frac{\bar{D}}{\hat{r}} \right)   \, + 1 \, \right)  
   \label{eqn.fl}    \\ &  \qquad        \quad \Rightarrow \quad \min_{k \leq \ell } \,   \frac{f(x_k') - f^*}{\hat{f} - f^*} \, \leq \, \epsilon \; .  \label{eqn.fm} \end{align} 

We prove Theorem~\ref{thm.aa}  by showing Algorithm A is ``equivalent'' to applying Algorithm 2 as above.  Likewise, we prove Theorem~\ref{thm.ab}  by showing algorithm B is equivalent to Algorithm 1. To be in position to establish these equivalences, however, we need a characterization of the supdifferentials $ \hat{\partial} \lambdamin(x,s,t) $, a characterization in terms of the original problem (\ref{eqn.fa}). \vspace{1mm}

\section{{\bf  A Practical Characterization of Supdifferentials, \\ and the Proofs of Theorems~\ref{thm.aa}  and \ref{thm.ab} }} \label{sect.g} 

We continue with the setting and notation developed in Section \ref{sect.f}.
\vspace{1mm} 

 We provide a practical characterization of the supgradients for the ``$ \lambda_{\min} $ function,'' a characterization expressed in terms of the original problem (\ref{eqn.fa}). The characterization provides the means for other algorithms designed using the framework to be implemented directly in terms of (\ref{eqn.fa}), avoiding computation of the conic reformulation. \vspace{1mm}

For the reader's convenience, we recall key sets appearing in Section~\ref{sect.a}  that also appear in the characterization below.  \vspace{1mm}

For $ x' \in \mathrm{bdy}(S) $, let
\[ G_1(x') := \left\{ \frac{-1}{ \lin v , \bar{e}  - x' \rin } \, v \mid  \, \vec{0} \neq  v \in N_S(x') \right\} \; . 
\]
For $ (x',t') \in \mathrm{\bdy}(\mathrm{epi}(f)) $, let         
\[ G_2(x',t') := \left\{ \frac{-1}{\lin v, \bar{e} - x' \rin + ( \hat{f} - t') \delta  } \, v \mid \,  (\vec{0},0) \neq (v,\delta) \in N_{\mathrm{epi}(f)}(x',t')   \right\} \; , \] 
where here, normality is with respect to the inner product that assigns pairs $ (x_1, t_1) $, $ (x_2, t_2) $ the value $ \lin x_1, x_2 \rin + t_1 t_2 $. 
\vspace{1mm}

Finally, letting $ \bdy_1 := \{ (x',t') \mid x' \in \bdy(S) \} $ and $ \bdy_2 := \bdy(\epi(f)) $, define for $ (x',t') \in \bdy_1 \cup \bdy_2 \; ,  $
\[ G(x',t') := \begin{cases}  G_1(x')  & \textrm{if $ (x',t') \in \bdy_1 \setminus \bdy_2 \; , $} \\
G_2(x',t') & \textrm{if $ (x',t') \in \bdy_2 \setminus \bdy_1 \; , $} \\
\textrm{hull}(G_1(x') \, \cup \, G_2(x',t') ) & \textrm{if $ (x',t') \in \bdy_1 \cap   \bdy_2 \; . $}   
\end{cases}    \]

\begin{prop}  \label{prop.ga} 
 Assume $ (x,1,t) \in \mathrm{Affine} $, and assume $ t < \hat{f} $. If $ (x',1,t') = \pi(x,1,t) $, then 
\[  P_{{\mathcal L}}( \, \hat{\partial}\lambda_{\min}(x,1,t) \, ) = \{ (P(v),0,0) \mid -v \in G(x',t') \} \; . \]
\end{prop}
\noindent {\bf Proof:}  Proposition \ref{prop.ca}  implies for all $ y \in {\mathcal E} \times \mathbb{R} \times \mathbb{R}  $,  
\begin{align}
   \hat{\partial} \lambda_{\min}(y) & = \{ \upsilon  \in - N_{{\mathcal K}}(y - \lambda_{\min}(y) e ) \mid \blin e, \upsilon  \brin = 1 \} \nonumber \\
  & = - \{ \upsilon \in N_{{\mathcal K}}(y - \lambda_{\min}(y) e ) \mid \blin e, \upsilon  \brin = - 1 \} \; . \label{eqn.ga} 
\end{align}

However, 
\begin{equation}  \label{eqn.gb} 
 \textrm{if $ y = (x,1,t) \in \mathrm{Affine}  $ and $ t < \hat{f} $,}
 \end{equation}
 then  $ \pi (y) $ is a positive multiple of $ y - \lambda_{\min}(y) e $ (using $ \lambda_{\min}(y) < 1 $ (Lemma~\ref{lem.bb}), and $ \pi(y) = e - \frac{1}{1 - \lambda_{\min}(y)}(y - e) $), and thus has the same normal cone to $ {\mathcal K} $. Hence, in this case, (\ref{eqn.ga})  is equivalent to
\begin{equation} \label{eqn.gc} 
  \hat{\partial} \lambda_{\min}(y) = - \{ \upsilon  \in N_{{\mathcal K}}(\, \pi (y) \,  ) \mid \blin e, \upsilon \brin = -1 \} \; . 
  \end{equation}
  
 As $ \mathrm{int}( {\mathcal K}) \neq \emptyset  $ and $ {\mathcal K} = {\mathcal K}_1 \cap {\mathcal K}_2 $, if $ y \in {\mathcal K} $ then
\[  N_{{\mathcal K}}(y) = \mathrm{hull}( \, N_{{\mathcal K}_1}(y) \cup N_{{\mathcal K}_2}(y) \, ) \; . 
\] 
   Thus, for $ y $ as in (\ref{eqn.gb}), using (\ref{eqn.gc})  we have
\begin{equation}  \label{eqn.gd} 
  \hat{\partial} \lambda_{\min}(y) = - \{ \upsilon  \in  \mathrm{hull}( \, N_{{\mathcal K}_1}(\pi(y)) \cup N_{{\mathcal K}_2}(\pi(y)) \, ) \mid \blin e, \upsilon \brin = -1 \} \; . 
  \end{equation}
\vspace{1mm}

Since $ e \in \mathrm{int}( {\mathcal K}_i) $ ($ i = 1,2 $), 
\[  y \in {\mathcal K}_i \, \wedge \, \upsilon \in N_{{\mathcal K}_i}(y) \, \wedge \,  \blin e, \upsilon \brin = 0 \quad \Rightarrow \quad  \upsilon = 0 \; . \]
Consequently, it follows from (\ref{eqn.gd})  that for $ y $ as in (\ref{eqn.gb}),
\[   
 \hat{\partial} \lambda_{\min}(y) =  -\mathrm{hull}( {\mathcal N}_1( \pi(y)) \cup {\mathcal N}_2(\pi (y)) ) \; , 
\] 
where for $ i = 1,2 $ and $ y' \in {\mathcal K}_i $, 
\[        {\mathcal N}_i(y') := \begin{cases} \{ v \in N_{{\mathcal K}_i}(y') \mid \blin e, \upsilon \brin = - 1 \} & \textrm{if $ y' \in \mathrm{bdy}({\mathcal K}_i$)}, \\ \{ \vec{0} \} & \textrm{if $ y' \in \mathrm{int}({\mathcal K}_i $)}.
\end{cases}  \]
Clearly, then,
\begin{equation} \label{eqn.ge} 
          P_{{\mathcal L}}( \hat{\partial} \lambda_{\min}(y)) =  - \mathrm{hull} \big( \,  P_{\mathcal L}(  {\mathcal N}_1( \pi(y))) \, \cup \, P_{{\mathcal L}}( {\mathcal N}_2( \pi(y))) \, \big) \; .
\end{equation}
              
Since $ S \times \mathbb{R} = \{ (x,t) \mid (x,1,t) \in {\mathcal K}_1 \} \; , $  if $ y' = (x',1,t') \in \bdy( {\mathcal K}_1) $ then 
\[ N_{{\mathcal K}_1}(y') = \{ (v, \gamma, 0) \mid v \in N_S(x') \, \wedge \lin x', v \rin  = 0 \} \; , \]
the equation to ensure perpendicularity to the ray through $ y' $.  Thus, a vector $ (v, \gamma, \delta) $ is an element of $ {\mathcal N}_1(y') $ if and only if $ v \in N_S(x') $, $ \delta = 0 $ and    
\[   \begin{array}{ccccc}
       \lin x', v \rin & + & \gamma  &  = & 0 \\
       \lin  \bar{e} , v \rin & + & \gamma  &  = & -1  
      \end{array} 
      \]
(note $ v $ cannot be $ \vec{0} $).       
However, for $ \vec{0} \neq v \in N_S(y') $, there is a unique scaling of $ v $ for which there exist $ \gamma $ and $ \delta $ satisfying the three equations, namely,
\[  v \mapsto \smfrac{-1}{\lin v, \bar{e} - x' \rin  }  (v, \,  - \lin x', v \rin , \,  0) \; . \]
The denominator in the scaling is negative, because $ \bar{e} \in \int(S) $.  It follows that  
\begin{equation}  \label{eqn.gf} 
  P_{\mathcal L} ( {\mathcal N}_1 ) = \{ \smfrac{-1}{\lin v, \bar{e} - x' \rin  } ( P(v), 0, 0) \mid \vec{0} \neq v \in N_S(x') \} \; . 
  \end{equation}

Similarly, for $ y' = (x',1,t') \in \bdy( {\mathcal K}_2) $, a vector $ (v,\gamma, \delta) $ is in $ {\mathcal N}_2(y') $ if and only $ ( v, \delta ) \in N_{ \epi(f)}( x', t') $ and 
\[   \begin{array}{ccccccc}
       \lin x', v \rin & + & \gamma  & + &  \delta t'  & = & 0 \\
       \lin \bar{e}, v  \rin & + & \gamma  & + & \delta \hat{f}  & = & -1 \
      \end{array} 
      \]
 (note $ (v, \delta) $ cannot be $ ( \vec{0},0) $).       
For $ ( \vec{0},0) \neq (v,\delta) \in N_{ \epi(f)}(x',t') $, the unique scaling for which there exists $ \gamma $ satisfying the equations is
\[ (v, \delta) \mapsto \smfrac{-1}{\lin v, \bar{e} - x' \rin + (\hat{f} - t') \delta}  (v, \, - \lin x', v \rin - t' \delta , \, \delta) \; .   \]
The denominator is negative, because $ (\bar{e}, \hat{f}) \in \int(\epi(f)) $.  It follows that
\begin{equation}  \label{eqn.gg} 
  P_{\mathcal L} ( {\mathcal N}_2 ) = \{ \smfrac{-1}{\lin v, \bar{e} - x' \rin + (\hat{f} - t') \delta} ( P(v), 0, 0) \mid (\vec{0},0) \neq (v, \delta) \in N_{\epi(f)}(x',t') \} \; . \end{equation}

Together, (\ref{eqn.ge}), (\ref{eqn.gf})  and (\ref{eqn.gg})  establish the proposition. \hfill $ \Box $
 \vspace{5mm}

\noindent 
{\bf Proofs of Theorems \ref{thm.aa}  and \ref{thm.ab}, and of (\ref{eqn.ag}):} We first establish Theorem~\ref{thm.aa}, and then remark on the few changes required to establish Theorem~\ref{thm.ab}  by the same approach. \vspace{1mm}

By ``Algorithm 2,'' we mean the application of Algorithm 2 to the conic optimization problem (\ref{eqn.fd}), as presented in the preceding section. \vspace{1mm}

Due to complexity bound (\ref{eqn.fi}), (\ref{eqn.fj})   established for Algorithm 2, it suffices to show the two algorithms are equivalent in that the iterates generated by Algorithm 2 are ``identical'' to those generated by Algorithm A. \vspace{1mm}

The initial iterates for Algorithm 2 are $ (x_0, 1, t_0) = (x_0',1,t_0') = (\bar{e}, 1, f( \bar{e}) $, and for Algorithm A are $ (x_0, t_0) = (x_0',t_0') = (\bar{e}, f(\bar{e})) $.  For beginning an inductive proof showing the two algorithms are equivalent, observe that the initial iterate $ (x_0, t_0) $ (resp., $ (x_0', t_0') $) for Algorithm A is obtained simply by eliminating the ``1'' from the initial iterate $ (x_0,1,t_0) $ (resp., $ (x_0',1,t_0') $) for Algorithm 2. \vspace{1mm}
 
For the inductive step, assume eliminating ``1'' from the iterate $ (x_k,1,t_k) $ (resp., $ (x_k',1,t_k') $) computed by Algorithm 2 results in the iterate $ (x_k,t_k) $ (resp., $ (x_k',t_k') $) computed by Algorithm A. 
\vspace{1mm}

The first iterate computed by Algorithm 2 in Step 1 is
    \[ (\tilde{x}_{k+1}, 1, t_k) =  (x_k, 1, t_k) + \smfrac{\epsilon }{2 \| P_{\mathcal L} \bar{g}  \|^2} P_{{\mathcal L}} \bar{g}  \quad \textrm{where } \bar{g}  \in \hat{\partial} \lambdamin(x_k', 1, t_k') \; , \]
whereas Algorithm A computes  
\[   \tilde{x}_{k+1} = x_k - \smfrac{\epsilon }{2 \| P g \|^2} Pg  \quad \textrm{where } g \in G(x_k', t_k')  \]
(and thereafter relies on the pair $ ( \tilde{x}_{k+1}, t_k) $ in the same manner that Algorithm 2 relies on the triple $ (\tilde{x}_{k+1}, 1, t_k) $).  Thus, due to the identity provided by  Proposition~\ref{prop.ga}, the algorithms can be considered equivalent in this computation. 
\vspace{1mm}

Algorithm 2 next computes $ (x_{k+1}', 1, t_{k+1}') := \pi(\tilde{x}_{k+1},1,  t_k) $, the point in the boundary of $ {\mathcal K} $ encountered when moving from $ (\bar{e}, 1, \hat{f}) $ in direction $ (\tilde{x}_{k+1},1 , t_k) - (\bar{e},1,\hat{f}) $. Due to the definition of $ {\mathcal K} $, however, this is the (only) feasible point $ (x,1,t) $  for which there exists $ \alpha \geq 0 $ such that
\[ x = x( \alpha) := \bar{e} + \alpha \cdot ( \tilde{x}_{k+1} - \bar{e}) \quad \textrm{and} \quad t = t( \alpha) := \hat{f} + \alpha \cdot ( t_k - \hat{f}) \]
and either $ x \in \mathrm{bdy}(S) $ or $ (x,t) \in \bdy(\epi(f))   $.  Thus, for the iterate $ (x_{k+1}', 1, t_{k+1}') $ computed by Algorithm 2, we have $ (x_{k+1}', t_{k+1}') = \pi'( \tilde{x}_{k+1}, t_k) $, where $ \pi' $ is defined by (\ref{eqn.ac}).  Thus, eliminating ``1'' from the iterate  $ (x_{k+1}', 1, t_{k+1}' ) $ computed by Algorithm 2 gives the iterate $ (x_{k+1}', t_{k+1}') $ computed by Algorithm A.  \vspace{1mm}

Observe, moreover, for Algorithm A, $ t_{k+1}' = \hat{f} + \alpha_{k+1} \cdot (t_k - \hat{f}) $, and thus,  
\begin{equation}  \label{eqn.gh} 
   \alpha_{k+1} = \frac{\hat{f} - t_{k+1}'}{\hat{f} - t_k} \; .
\end{equation}

Algorithm A decides how to define its iterate $ (x_{k+1}, t_{k+1}) $ based on whether $ \alpha_{k+1} \geq 4/3 $, whereas Algorithm 2 decides how to define its iterate $ (x_{k+1},1,t_{k+1}) $ based on whether $ \hat{f} - t_{k+1}' \geq \smfrac{4}{3} ( \hat{f} - t_k) $.  Thus, due to (\ref{eqn.gh})  and the equivalence of the two algorithms up until now, it is simple to verify that removing ``1'' from the iterate $ (x_{k+1},1,t_{k+1}) $ for Algorithm 2 gives the iterate $ (x_{k+1}, t_{k+1}) $ for Algorithm A. \vspace{1mm}

Similarly, the proof of Theorem~\ref{thm.ab}  is accomplished by showing Algorithm B is ``equivalent'' to Algorithm 1, and making use of the complexity result (\ref{eqn.fk}), (\ref{eqn.fl}), (\ref{eqn.fm})    for Algorithm 1.  The proof of equivalence is virtually identical to the one above, the main difference being that the scalar $ \alpha_k $ used by Algorithm B has to be related to the analogous scalar used by Algorithm 1 -- in particular, for equivalence, the identity $ \frac{\alpha(x_k, f^*) - 1}{\alpha(x_k, f^*)} = \lambda_{\min}(x_k,1,f^*) $ is needed. \vspace{1mm}

The identity, however, is easily established.  Indeed, $ \alpha(x_k,f^*) $ is the smallest scalar $ \alpha $ for which the point $ (x,t) = ( \bar{e}, \hat{f}) + \alpha \cdot ( (x_k, f^*) - (\bar{e}, \hat{f})) $  satisfies either $ x \in \mathrm{bdy}(S) $ or $ (x,t) \in \bdy(\epi(f))  $. Due to the definition of $ {\mathcal K} $, however, this  is readily seen to equal the smallest scalar $ \alpha $ for which $ e + \alpha \cdot ((x_k,1,f^*) - e) $ lies in the boundary of $ {\mathcal K} $.  Hence, from (\ref{eqn.bda}), $ \alpha(x_k,f^*) = \frac{1}{1 - \lambda_{\min}(x_k,1,f^*)} $. Rearrangement establishes the desired identity. \vspace{1mm}

Finally, immediately after the statement of Theorem~\ref{thm.ab}, it is claimed in (\ref{eqn.ag})  that Algorithm B converges linearly when $ f $ is a piecewise linear function and $ \mathrm{Feas} $ is polyhedral.  This is immediate from the fact that Algorithm 1 converges linearly when $ {\mathcal K} $ is a polyhedral cone (Corollary~\ref{cor.db}).  
  \hfill $ \Box $
  \vspace{3mm}

{\bf Acknowledgements:}  The author expresses gratitude to the reviewers and the associate editor, for careful attention to detail and for suggestions on how to position the paper.  The author thanks Yurii Nesterov for encouragement at an early, critical stage of the research, Rob Freund for discussions and feedback influencing the overall manner of presentation, and Wei Qian for conversations affecting the deduction and presentation of results for general convex optimization.

\bibliographystyle{plain}
\bibliography{efficient_and_general_subgradient-methods}

  \end{document}